\documentclass[a4paper]{amsart}

\RequirePackage{amsmath, amssymb, amsthm, amsfonts}
\usepackage[all]{xy}
\usepackage{enumitem}
\usepackage{hyperref}
\usepackage{fullpage}
\usepackage{hyperref}
\usepackage[mathscr]{eucal}

\newif\ifpdf
\ifx\pdfoutput\undefined \pdffalse \else \ifnum\pdfoutput=0
\pdffalse \else \pdfoutput=1 \pdftrue \fi \fi

         \ifpdf
         \usepackage[pdftex]{graphicx}
         \else
         \usepackage{graphicx}
         \fi
\newtheorem{theo}{Theorem}[section]
\newtheorem{lem}[theo]{Lemma}
\newtheorem{prop}[theo]{Proposition}
\newtheorem{cor}[theo]{Corollary}
\newtheorem{defin}[theo]{Definition}
\theoremstyle{definition}

\newtheorem{example}[theo]{Example}

\newtheorem{rem}[theo]{Remark}

\newcommand{\N}{\mathbb{N}}

\newcommand{\Z}{\mathbb{Z}}
\newcommand{\Q}{\mathbb{Q}}

\newcommand{\C}{\mathbb{C}}

\newcommand{\PP}{\mathbb{P}}

\newcommand{\Aut}{\rm Aut}
\newcommand{\id}{\rm id}

\newcommand{\rk}{\mathrm{rk}}

\newcommand{\ZZ}{\mathbb{Z}}

\newcommand{\CC}{\mathbb{C}}

\newcommand{\lr}{\longrightarrow}

\newcommand{\ra}{\rightarrow}

\title{K3 surfaces with a non-symplectic automorphism and product-quotient surfaces with cyclic groups}

\author{Alice Garbagnati and Matteo Penegini}

\begin{document}

\subjclass[2010]{Primary 14J28; Secondary  14J10, 14J50, 14D06.}
\keywords{K3 surfaces, automorphisms of K3 surfaces, product-quotient surfaces.\\
Both authors are partially supported by PRIN 2010--2011 ``Geometria delle variet\`a algebriche"}

\maketitle

\begin{abstract}
We classify all the K3 surfaces which are minimal models of the
quotient of the product of two curves $C_1\times C_2$ by the
diagonal action of either the group $\Z/p\Z$ or the group
$\Z/2p\Z$. These K3 surfaces admit a non-symplectic automorphism
of order $p$ induced by an automorphism of one of the curves
$C_1$ or $C_2$. We prove that most of the K3 surfaces
admitting a non-symplectic automorphism of order $p$  (and in fact
a maximal irreducible component of the moduli space of K3 surfaces
with a non-symplectic automorphism
of order $p$) are obtained in this way.\\
In addition, we show that one can obtain the same set of K3
surfaces under more restrictive assumptions namely
one of the two curves, say $C_2$, is isomorphic to a rigid
hyperelliptic curve with an automorphism $\delta_p$ of
 order $p$ and the automorphism of the K3 surface is induced by $\delta_p$.\\
Finally, we describe the variation of the Hodge structures of the surfaces constructed and 
we give an equation for some of them.
\end{abstract}

\section{Introduction}\label{sec.intro}

One of the main themes of interest in the study of K3 surfaces $S$
regards their automorphisms. We call an automorphism ${\bf g}$ of $S$ \emph{non-symplectic} if it acts
non-trivially on the nowhere vanishing holomorphic 2-form
$\omega$. In the case $|{\bf g}|=p$, a prime number,  ${\bf g}(\omega)=\zeta_p \omega$, where $\zeta_p$
is a primitive $p$-th root of unity.
The pairs $(S, {\bf g})$ are
quite rare, in the sense that there are strict restrictions on
both $p$, which must be smaller then or equal to $19$, and on the
K3 surface $S$, which cannot be generic in the moduli space. More
precisely, the families of K3 surfaces with a non-symplectic automorphism of odd prime order have a
finite number of connected components and the biggest of them has
dimension $9$ while the moduli space of the K3 surfaces is
20-dimensional.

Several authors worked on the classification of K3 surfaces
admitting a non-symplectic automorphism of odd prime order (see,
e.g., \cite{AS08}, \cite{AST},\cite{Kondo}, \cite{OZ1},
\cite{OZ2}, \cite{OZ3}) and eventually their complete
classification is given in \cite{AST}. The classification is based on the following procedure: First, a non-symplectic automorphism {\bf g} of order $p$ acting on a K3 surface $S$, determines an action on the lattice $T:=(H^2(S,\Z)^{\bf g})^{\perp}) \subset H^2(S,\Z)$ . The lattice $T$ satisfies several conditions. Second, one lists all the lattices with such conditions. Third,  for each lattice $T$ in the list one has to actually construct an example of a K3 surface admitting a non-symplectic automorphism ${\bf g}$ and such that $T\simeq (H^2(S,\Z)^{\bf g})^{\perp})$. Each example is given by an {\it ad hoc} construction.

The aim of this paper is to give a systematical way to construct
most of these K3 surfaces by showing that they are the minimal models of
product-quotient surfaces (i.e., of the minimal resolution of the
quotient $(C_1\times C_2)/G$ where the $C_i$ are curves of genus $g(C)\geq 1$ and $G$ is a
finite group acting diagonally on the product, see Definition
\ref{defin PQ}). In addition, as we have already observed, the
pairs $(S, {\bf g})$ are quite "special", we give here a
geometrical interpretation of this "speciality" in many cases:  $S$ is the minimal model of a product-quotient surface. Moreover, we prove that there is a curve which seems to play a central r\^{o}le in this construction. We define
the curve $D_p$ as the hyperelliptic curve with equation
$v^2=u^p-1$. It clearly admits an automorphism $\delta_p$ of order $p$, acting
on $u$ as the multiplication by $\zeta_p$, and an automorphim
$\tau_p$ of order $2p$ which is the composition of $\delta_p$ with
the hyperelliptic involution. We prove the following theorem.

\begin{theo}\label{theo main} If $S$ is a K3 surface admitting a non-symplectic automorphism of order $p=3$
(resp. $p=5,7,11,13$, $p=17,19$) whose fixed locus contains at
least 2 (resp. 1, 0) curves, then it is the minimal model of a
resolution of the quotient $(C_1\times D_p)/({\bf g_1}\times
\tau_p)$, where ${\bf g_1}$ is an automorphism of $C_1$ of order
2p. The non-symplectic automorphism of order $p$ on $S$ is induced by the automorphism $\id\times \delta_p$.
\end{theo}

Any K3 surface $S$ constructed in the Theorem \ref{theo main} admits two isotrivial families of curves, whose general member is isomorphic to $C_1$ and $D_p$ respectively. The non-symplectic automorphism of order $p$ on $S$ is in fact given by the action of an automorphism on each member of one of these families.

We observe that the K3 surfaces given by the Theorem \ref{theo main} admit a
non-symplectic automorphism of order $2p$ too, the one induced by
$\id\times \tau_p$. This gives a geometrical explanation of the
following significative result on non-symplectic automorphisms: If
a K3 surface $S$ admits a non-symplectic automorphism of order
$p$, under certain conditions on the fixed locus it admits in fact
a purely non-symplectic automorphism of order $2p$ (see
\cite[Theorems 1.6, 1.7]{GS} and \cite{Dillies} for the precise
statement).

The proof of the Theorem \ref{theo main} is based on the construction of the K3 surfaces that we now explain. We first bound
the genus of ther curves $C$ with a cyclic group of automorphisms $G$
of order $p$ (resp. $2p$), having the properties that $C/G \simeq \PP^1$,
and such that there exists an eigenspace $H^{1,0}(C)_{\zeta^i_p}$ (resp. $H^{1,0}(C)_{-\zeta^i_p}$) of
dimension $1$ of the induced action. This is achieved by
exploiting the Riemann Existence Theorem, the Holomorphic
Lefschetz Fixed-Point formula, and the Chevalley--Weil formula.
Second, we classify all these curves and we couple them by choosing the action of $G$ in such a way that the singular surfaces $(C_1 \times C_2)/G$ have $p_g=1$ and
$q=0$. Third, we resolve the singularities and we get
product-quotient surfaces $X$ which are not minimal models, but in two cases. We observe that $K^2_X$ could be very negative. After having found all the $(-1)$-curve
on $X$ (this is a quite delicate task see e.g., \cite{BP11}) we
carefully contract them to produce a minimal model $S$ of $(C_1\times C_2)/G$. Finally, we
prove that $S$ is a K3 surface.

As a byproduct of the proof we classify the K3 surfaces
which are minimal models of the product-quotient surfaces with the
groups $\Z/p\Z$ and $\Z/2p\Z$. Hence, we obtain also a "negative" result: if a K3
surface does not satisfy the hypothesis of Theorem \ref{theo main}
(i.e.,\ either it does not admit a non-symplectic automorphism of
order $p$, or it admits a non-symplectic automorphism of order
$p$, but its fixed locus does not satisfy the condition of the
Theorem \ref{theo main}), then this K3 surface is not the minimal model of a
product-quotient with group $\Z/2p\Z$. However, we cannot exclude
that such a K3 surface is the minimal model of a product-quotient
with a different group. Indeed, we are aware that this is the case for at least certain families of K3 surfaces admitting a
non-symplectic automorphism of order 3. We shall analyze this problem in
a forthcoming article. \\

This paper is organized as follows.

In Section \ref{sec.curves} we briefly recall three classical
results: the Riemann Existence Theorem, the Holomorphic Lefschetz
Fixed-Point formula, and the Chevalley--Weil formula. We establish
the upper and lower bound for the genus of a curve $C$ with a
cyclic group of automorphisms of odd prime order, and the property
that there exists an eigenspace $ H^{1,0}(C)_{\zeta^i_p}$ of
dimension $1$. We state similar results for
cyclic groups of automorphisms of order $2p$. Moreover, we give some explanatory
examples introducing the curve $D_p$.

Section \ref{Sec Surfaces} is divided into three parts. In the
first part we give the definition of product-quotient surfaces and
we recall the properties of these surfaces that are needed for our
purposes. In the second part we calculate the Hodge numbers of
product-quotient surfaces. In the last subsection we describe the
automorphisms of the minimal model of a product-quotient surface.

In Section \ref{Sec K3 sur} we describe, first, the procedure we
used to construct product-quotient surfaces with group either
$\ZZ/p\ZZ$ or $\ZZ/2p\ZZ$ and $p_g=1$ and $q=0$. Second, we give a
method to prove that these surfaces are K3 surfaces.

Sections \ref{sec: PQ with Zp} and \ref{sec: PQ with Z2p} present the main results of the paper on K3 surfaces which are minimal models of product-quotient surfaces with group either $\Z/p\Z$ or $\Z/2p\Z$ respectively. Moreover, in
these two sections one can find the tables with the surfaces we constructed (see \verb|Table 1| and \verb|Table 2|).

In Section \ref{sec: equation} we give the equations for the singular models of the K3 surfaces
constructed for $p=3,5,7,11$. In particular if $p=5,7,11$ we describe these surfaces as hypersurfaces in weighted projective spaces.

In the last section we describe the variation of the Hodge structures of the K3 surfaces constructed in terms of the Hodge structure of $H^1(C_1)$ relating them by a half twist.

In the Appendix there is the \verb|MAGMA| script of the program that we used. There are essentially two programs \verb|Surfacesp| and \verb|t1t2PtsSurfaces|. The former one gives a list of all product-quotient surfaces $X$ with group $\ZZ/p\ZZ$ ($p$ any odd prime), $p_g(X)=1$ and $q(X)=0$, as well as the singularities of $(C_1\times C_2)/({\bf g_1}\times {\bf g_2})$, this program becomes very slow as $p$ increases. The latter program gives a similar list, it is much faster, it works also for the group $\ZZ/2p\ZZ$, but it requires two additional data which are the number of ramification points of the two coverings $C_i \ra C_i/G$.
\bigskip

{\bf Notation:} We work over the field of complex numbers $\CC$.
We will denote by $\zeta_n:=e^{\frac{2 \pi i}{n}}$ a primitive $n$-th root of unity.\\
By ``curve" or ``Riemann surface"  we mean a projective, non-singular curve $C$, we denote by $H^{1,0}(C)=H^0(C,\Omega^1_C)$ and by $g(C):=h^{1,0}(C):=dimH^{1,0}(C)$ the \emph{genus} of the curve.

By ``surface'' we mean a projective, non-singular surface $S$, and
for such a surface $\omega_S=\mathcal{O}_S(K_S)$ denotes the canonical
class, $p_g(S)=h^{2,0}(S)=h^0(S, \, \omega_S)$ is the \emph{geometric genus},
$q(S)=h^{1,0}(S)=h^1(S, \, \omega_S)$ is the \emph{irregularity} and
$\chi(S)=1-q(S)+p_g(S)$ is the \emph{Euler-Poincar\'e
characteristic}. The \emph{Noether formula} is $12 \chi(S)=K^2_S+e(S)$, where $e(S)$ is the Euler number of $S$.

By abuse of notation by ``$(-1)$-curve'' we mean a curve $C$ with $C \simeq \PP^1$ and $C^2=-1$.\\

\textbf{Acknowledgments}: We warmly thanks Bert van Geemen and Roberto Pignatelli for several essential suggestions and useful discussions.


\section{Curves}\label{sec.curves}


This section is devoted to recall some classical results: a reformulation of the Riemann Existence Theorem, the Holomorphic Lefschetz Fixed-Point formula, and the Chevalley--Weil formula, as well as fixing the notation. Moreover, here we establish all the properties for the curves that we will use to construct K3 surfaces in the next sections.

\subsection{The Riemann Existence Theorem}

\begin{defin}\label{defin.orbifold.grp}
Let $m_1, \dots , m_r$ be positive integers with $m_i \geq 2$ for
all $i$. A \emph{polygonal group} of type $(m_1, \dots , m_r)$ is
a group presented as follows:
\begin{equation*}  \Gamma(m_1, \dots , m_r):=\langle\gamma_{1},
\ldots , \gamma_{r} |  \gamma^{m_1}_{1}=\dots=
\gamma^{m_r}_{r}=\gamma_{1} \cdot \ldots \cdot \gamma_{r} =1
\rangle.
\end{equation*}
\end{defin}

\begin{defin}\label{defin.adm.epi}
Let $\Gamma$ be a polygonal group and $G$ be a finite group. An
epimorphism $\theta: \Gamma(m_1,...,m_r) \rightarrow G$ is called
\emph{admissible} if $\theta (\gamma_i)=x_i$ has order exactly $m_i$ for all
$i$.
If an admissible epimorphism exists, then
the image $(x_1, \ldots ,x_r)$ of a set of generators of $\Gamma$ is called a \emph{spherical system of generators} for $G$.
\end{defin}

The following is a reformulation of {\em Riemann's Existence
Theorem} (see e.g., \cite{M95} chapter III, sections 3 and 4):

\begin{theo}\label{theo: Riemann} A finite group $G$ acts as a group of automorphisms on some
compact Riemann surface $C$ such that $C/G \simeq \mathbb{P}^1$ if and only if there are natural numbers $m_1,
\ldots , m_r$, and an admissible epimorphism
$$\theta \colon
\Gamma(m_1,\ldots, m_r) \rightarrow G.$$ The genus $g(C)$ is determined by the Riemann--Hurwitz relation:
\begin{equation}\label{form.RH} 2g(C) - 2 = |G|\left(-2 +
\sum_{i=1}^r \left(1 - \frac{1}{m_i}\right)\right).
\end{equation}
\end{theo}
The $G-$cover $C \rightarrow \mathbb{P}^1$ is
branched in $r$ points $p_1,\dots,p_r$ with branching number
$m_1,\dots,m_r$, respectively. Moreover, the cyclic subgroups
$\left\langle x_{i}\right\rangle$ and their conjugates are the
non-trivial stabilizers of the action of $G$ on $C$.

In this paper we always assume $G$ to be cyclic, so any spherical system of generators is of the form $(x^{\xi_1}, \ldots ,x^{\xi_r})$, where $\xi_i\in \Z/n\Z$.
Moreover, we are interested only in the cases:
\[
G = \ZZ/n\ZZ\mbox{ with }n=p\mbox{ or }n=2p\mbox{ and }p\mbox{ prime.} 
\]
From now on we shall assume that we are in one of these cases. Notice that some of the results are true also for more general cases.

As an automorphism group of $C$, $G=\langle {\bf g} \rangle $
gives a conformal self-mapping ${\bf g}\colon C \rightarrow C$ of
order $n$. 
Please notice
that the abstract group $G$ has also other realizations: as image of an admissible epimorphism $\langle x \rangle$ or as
local action $\langle \zeta_n \rangle$ near a point $P \in C$.
Suppose that ${\bf g}^k$ fixes a point $P \in C$, then in a suitable
local coordinate $z$ near $P$ we must have ${\bf g}^k(z)=\zeta_n^i z$. Thus
${\bf g}^k$ is locally a rotation at $P$ and the rotation angle is determined in the following proposition.

\begin{prop}\cite[Theorem 7]{Ha71}\label{prop: rotation consts} Let $C$ be a curve associated to the spherical system of generators $(x^{\xi_1},\ldots, x^{\xi_r})$. Let $P_j$ be a point with a non trivial stabilizer in $\langle {\bf g} \rangle \subset {\rm Aut}(C)$ of order $m$ generated by ${\bf g}^{\xi_j}$. Then ${\bf g}^{\xi_j}(z)=\zeta_n^{\eta_j}z$ where $z$ is the local coordinate near the point $P_j$ and $\xi_j\eta_j\equiv n/m  \mod n$ and $0<\eta_j < m$.
\end{prop}
Assume that $n=p$ with $p$ prime, denote by $a_i$ the number of
ramification points of $C$ where the local action is given by
$\zeta_p^i$. Notice that in this case we have total ramification
over each branch point. The condition of Theorem \ref{theo: Riemann} is equivalent to the condition
\begin{equation}\label{eq.monodromy.cond}
\sum^{p-1}_{i=i}i^{-1}a_i \equiv 0 \textrm{ mod }p,
\end{equation}
where the inverse of $i$ is taken$\mod p$. We shall call this condition the \emph{numerical monodromy condition}.

\subsection{The Holomorphic Lefschetz Fixed-Point Formula}

Let $\ZZ/n\ZZ \simeq G \subset Aut(C)$, then there is an induced linear action of $G$ on the Dolbeault cohomology $H^{*,*}(C)$. Let us recall the holomorphic analogous of the Lefschetz Fixed-Point formula.

\begin{theo}{\rm(Holomorphic Lefschetz Fixed-Point Formula \cite[p. 426]{GH})} Let $G$ be a group of automorphism of a smooth complex curve $C$. Let ${\bf g} \in G$ then it holds
\begin{equation*}
1- Tr ({\bf g} \mid_{H^{0}(X,\Omega^1_X)})=\sum_{\begin{array}{c} \scriptstyle P \in C \\
\scriptstyle {\bf g}(P)=P  \\
 \end{array}} \frac{1}{1-\zeta_P},
\end{equation*}
where $\zeta_P$ is the local action of ${\bf g}$ in $P$.
\end{theo}

In particular in the case $G \simeq \ZZ/p\ZZ$ we have a decomposition in eigenspaces
\begin{equation*}
H^{1,0}(C) = \bigoplus^{p-1}_{j=1} H^{1,0}(C)_{\zeta^j_p}.
\end{equation*}
The dimension of each eigenspace will be denoted by
\begin{equation*}
\alpha_j:=dim H^{1,0}(C)_{\zeta^j_p}.
\end{equation*}

With this piece of notation we can rewrite the Holomorphic Lefschetz Fixed-Point Formula in the following way.

\begin{cor} Let $C$ be a curve and $G= \langle \zeta_p \rangle \simeq \ZZ/p\ZZ$ then it holds

\begin{equation}\label{eq.lef_curve_gen}
-\sum^{p-1}_{j=1}\zeta^j \alpha_j=\sum^{p-1}_{l=1}a_l\frac{1}{1-\zeta^l}-1.
\end{equation}
\end{cor}

\begin{prop}{\rm(The Chevalley--Weil Formula)}\label{prop: alpha ai} Let $C$ be a cyclic $p:1$ cover $\mathbb{P}^1$. Let $a_i$ be the number of ramification points of $C$ where the local action is
given by $z \mapsto \zeta_p^iz$, and $\alpha_j:=dim H^{1,0}(C)_{\zeta^j_p}$. Then  $\alpha_j$ are determined by the $a_i$'s:
\begin{equation}\label{eq.rel.a.alpha}
 \alpha_{p-r}=-1+\frac{1}{p}\sum^{p-1}_{l=1}a_l l^{-1}(p-r),
\end{equation}
where the product $l^{-1}(p-1)$ is taken$\mod p$, and $r=1, \ldots ,n-1$.
\end{prop}
\begin{proof} First let us consider the LHS of \eqref{eq.lef_curve_gen}
\begin{equation*}
-\sum^{p-1}_{j=1}\zeta^j \alpha_j= -\sum^{p-2}_{j=1}\alpha_j\zeta^j+\alpha_{p-1}+\sum^{p-2}_{j=1}\alpha_{p-1}\zeta^j=
\alpha_{p-1}+\sum^{p-2}_{j=1}(\alpha_{p-1}-\alpha_j)\zeta^j.
\end{equation*}
We can rewrite the above equation as
\begin{equation}\label{eq.coeff.zeta.alpha}
\left[ \begin{array}{cccccc}
-1 &  0 & \ldots & \ldots & 1 \\
0  & -1 &  0     & \ldots & 1 \\
\vdots  &   &   \ddots  &  & 1 \\
0  &   &    & -1 & 1 \\
0 & \ldots & \ldots & 0  & 1 \\
\end{array} \right] \left[
\begin{array}{c} \alpha_1 \\
\alpha_2 \\
\vdots \\
\vdots \\
\alpha_{p-1}
\end{array} \right] =
\left[
\begin{array}{c} \textrm{coeff. of } \zeta \\
\textrm{coeff. of } \zeta^2 \\
\vdots \\
\vdots \\
\textrm{coeff. of } 1
\end{array} \right].
\end{equation}
Second looking at the RHS of \eqref{eq.lef_curve_gen} we have,
\begin{equation}\label{eq.coef.zeta.a}
\sum^{p-1}_{l=1}a_l\frac{1}{1-\zeta^l}=\frac{1}{p}\sum^{p-1}_{l=1}a_l\big(\sum^{p-1}_{k=0}(p-1-k)\zeta^{lk}\big)=\frac{1}{p}\sum^{p-1}_{l=1}a_l\big(\sum^{p-1}_{h=0}(p-1-hl^{-1})\zeta^{h}\big)=
\end{equation}
\begin{equation*}
=\frac{1}{p}\sum^{p-1}_{l=1}a_l\big[\sum^{p-2}_{h=0}(p-1-hl^{-1})\zeta^{h}-\big(p-1-(p-1)l^{-1}\big)\sum^{p-2}_{h=0}\zeta^{h}\big]=\frac{1}{p}\sum^{p-1}_{l=1}a_l\big(\sum^{p-2}_{h=0}l^{-1}(p-h-1)\zeta^{h}\huge),
\end{equation*}
where the product $l^{-1}(p-h-1)$ is taken$\mod p$.
Comparing the coefficient of $\zeta^0$ from the LHS and the RHS of \eqref{eq.lef_curve_gen} we obtain
\begin{equation}\label{eq: alpha n-1}
\alpha_{p-1}=-1+\frac{1}{p}\sum^{p-1}_{l=1}a_ll^{-1}(p-1).
\end{equation}
More generally we have by \eqref{eq.coeff.zeta.alpha}
\begin{equation*}
\alpha_{p-r}=\alpha_{p-1}-\textrm{coeff of }\zeta^r.
\end{equation*}
By substituting the expression of the coefficient of $\zeta^r$ given in \eqref{eq.coef.zeta.a} we obtain \eqref{eq.rel.a.alpha}.
\end{proof}

Notice that the above proposition is a special case of the well known Chevalley--Weil Formula, see \cite{CW}.

The above proposition has many important applications.
\begin{cor}\label{cor: genus ai}\label{cor: max genus} Let $r$ be the  number of ramification points for $\pi\colon C \rightarrow C/(\ZZ/p\ZZ) \simeq \PP^1$.
If there exists $i\in\{1,\ldots, p-1\}$ such that $\alpha_i=1$, then $(p-1)/2\leq g(C)\leq (p-1)^2$ and $3\leq r\leq 2p$.
\end{cor}
\begin{proof} By \eqref{eq.rel.a.alpha} we have
\begin{equation*}
g(C)=\sum^{p-1}_{k=1} \alpha_k=1-p+\sum^{p-1}_{k=1}\big[\frac{1}{p}\sum^{p-1}_{l=1}a_ll^{-1}(p-k)\big].\\
\end{equation*}
Notice that $\sum^{p-1}_{k=1}l^{-1}(p-k)$ is the sum of the first $p-1$ integers, hence it is equal to $p(p-1)/2$, and this gives $g(C)$.
The equation $\alpha_i=1$ gives conditions on the $a_i$'s. Substituting these values in $g(C)$ one obtains $(p-1)^2-\sum_i\lambda_ia_i$ with $\lambda_i\in\N$, which is less then or equal to $(p-1)^2$. The minimal genus is realized by the curve in Example \ref{ex: curve Dp}.\\
The case $p=3$ was already studied in \cite{vG}.
\end{proof}

Notice that by \eqref{eq: alpha n-1} the condition: $\alpha_{p-1}$ being a natural number, is equivalent to the numerical monodromy condition \eqref{eq.monodromy.cond}. Hence, by the Riemann Existence Theorem we have a curve with group of automorphism $\ZZ/p\ZZ$ and quotient $\PP^1$ once we provide a integral soltution of \eqref{eq: alpha n-1}.

\begin{example}\label{ex: curve Dp} Let $G\simeq \Z/p\Z$. The triple $(x^{p-1},x^{p-1},x^2)$ is a spherical system of generators for $G$. Hence, by Theorem \ref{theo: Riemann} there exists a curve, $D_p$, such that $D_p$ is a $p:1$ cover of $\mathbb{P}^1$ branched in 3 points and the cover automorphism $\delta_p$ acts locally as $\zeta_p^{p-1}$ near two fixed points and as $\zeta_p^{(p+1)/2}$ near the other one. This means that the only non zero $a_i$'s are $a_{p-1}=2$, $a_{(p+1)/2}=1$. For every choice of three points in $\mathbb{P}^1$ there exists an involution of $\mathbb{P}^1$ switching the first two points and fixing the other, so there exists an involution of $\mathbb{P}^1$ acting in this way on the branch points of $D_p\ra \mathbb{P}^1$. This induces an involution on $D_p$, which will be denoted by $\iota_p$. We observe that $\iota_p\delta_p=\delta_p\iota_p$ and we denote by $\tau_p=\iota_p\delta_p$. An equation of $D_p$ and the corresponding equation for its automorphisms are: \begin{equation*}\label{formula Dp}D_p\colon u^p=v^2-1, \quad \delta_p:(u,v)\mapsto (\zeta_pu,v),\ \iota_p:(u,v)\mapsto (u,-v),\ \tau_p:(u,v)\mapsto (\zeta_pu,-v).\end{equation*}
The genus of $D_p$, computed by the Riemann--Hurwitz formula, is $g(D_p)=(p-1)/2$. The curve $D_p$ is hyperelliptic over $\mathbb{P}^1_{[u]}$ and so a basis for $H^{0,1}(D_p)$ is given by $\{u^jdu/v\}$, $j=0,\ldots , (p-3)/2$. Therefore the eigenspaces decomposition of $H^{0,1}(D_p)$ for the induced action of $\delta_p$ is $\alpha_i=1$ if $i=1,\ldots, (p-1)/2$ and $\alpha_i=0$ if $i=(p+1)/2,\ldots, p-1$. This agrees with the results stated in Proposition \ref{prop: alpha ai}.\\
Finally we observe that $D_3$ is the elliptic curve with complex multiplication of order 3, which is associated to the automorphism $\delta_3$. \end{example}

There are analogous results for $G\simeq \ZZ/2p\ZZ$. 
The calculations are similar to the ones done above, which is why we omit it here. To get the formulae one can use the  \verb|MAGMA| script \verb|MaxGenus| in the Appendix. Despite the similarities, there is a difference between the case $\ZZ/p\ZZ$ and the case $\ZZ/2p\ZZ$. In the latter case we succeeded in finding a maximal genus only if we ask that the eigenspace $H^{0,1}(C)_{\zeta^i_{2p}}=1$ is relative to an element of maximal order. Indeed, we obtain the following.

\begin{cor}\label{rem: max genus 2p}
Let $r$ be the  number of ramification points for $\pi\colon C \rightarrow C/(\ZZ/2p\ZZ) \simeq \PP^1$.
If there exists $i\in\{1,\ldots, p-1\}$ such that $dim(H^{1,0}(C)_{-\zeta_p^i})=1$, then $(p-1)/2\leq g(C)\leq (2p-1)^2$ and $3\leq r\leq 4p$.
\end{cor}

\begin{example}\label{ex: Dp autom 2p} Let us consider the curve $D_p$ and the automorphism $\tau_p$ of order $2p$ (cf. Example \ref{ex: curve Dp}) whose associated spherical system of generators is $(x^{p+2}, x^{2p-2}, x^p)$. The action of $\tau_p$ on the form $u^idu/v$ is the multiplication by $-\zeta_p^{i+1}$, $0\leq i\leq (p-3)/2$. Hence, the eigenspaces decomposition of $H^{1,0}(D_p)$ for $\tau_p$ is the following: if $0\leq i\leq (p-3)/2$ and $i\equiv 0 \mod 2$, then $\alpha_{i+1}=1$; if $0\leq i\leq (p-3)/2$ and $i\equiv 1 \mod 2$, then $\alpha_{p+i+1}=1$; $\alpha_j=0$ otherwise. We observe that $\tau_p$ switches two of the fixed points of $\delta_p$ and fixes the third fixed point of $\delta_p$. Clearly the two points switched by $\tau_p$ are the ones where the local action of $\delta_p$ is $\zeta_p^{p-1}$ and the point fixed both by $\tau_p$ and by $\delta_p$ is the one where the local action of $\delta_p$ is $\zeta_p^{(p+1)/2}$. The automorphism $\tau_p^2$ coincides with $\delta_p^2$ and this identifies the local action of $\tau_p$: $\tau_p$ fixes one point with local action $-\zeta_p^{(p+1)/2}$; $\tau_p^2$ fixes other two points (switched by $\tau_p$) with local action $\zeta_p^{p-2}$; $\tau_p^p$ fixes other $p$ points (permuted by $\tau_p$) with local action $-\zeta_p^p=-1$.
\end{example}


\section{Surfaces}\label{Sec Surfaces}
In this section we recall the properties of the product-quotient surfaces and their minimal models and we calculate their numerical invariants.
\subsection{Product-Quotient Surfaces.}


Let us consider two curves $C_1$ and $C_2$ of genera greater then or
equal to $1$, and their product $C_1\times C_2$. Then
$\Aut(C_1)\times \Aut (C_2)\subset \Aut (C_1\times C_2)$. One can
say even more in the case $g(C_i) \geq 2$.
\begin{lem}\cite[Corollary 3.9]{cat00} Let us assume $g(C_1)\geq 2$ and $g(C_2)\geq 2$. If $C_1\not\simeq C_2$, then $\Aut(C_1\times C_2)=\Aut(C_1)\times\Aut(C_2)$, otherwise $\Aut(C_1\times C_2)=\left(\Aut(C_1)\times\Aut(C_2)\right)\rtimes\ZZ/2\ZZ$
\end{lem}

Let $G \subset \Aut(C_1)\times \Aut (C_2)$ be a finite group and
consider $(C_1 \times C_2)/G$, where $G$ acts diagonally on the
product $C_1 \times C_2$.

\begin{defin}\label{defin PQ}
The minimal resolution $X$ of the singularities of $(C_1 \times
C_2)/G$, where $G$ is a finite group with a diagonal action on the
direct product of two smooth curves $C_1$ and $C_2$ of respective
genera at least $1$, is called a \emph{product-quotient surface
with group $G$}.

We call $(C_1 \times C_2)/G$ the \emph{quotient model} of the product-quotient surface.
\end{defin}

Product-quotient surfaces were intensively studied, we refer to
\cite{BP11}, \cite{MP} and \cite{BCGP}  for a detailed account of them. We recall only the
facts that are important for our purposes.

\begin{rem}\label{rem: PQ properties} The following facts hold:
\begin{enumerate}
\item  There are only finitely many points on $C_1 \times C_2$ with non trivial stabilizer, which is cyclic. Therefore the quotient model has only a finite number of cyclic quotient singularities.

A \emph{cyclic quotient singularity} is locally analytic isomorphic to the quotient of $\CC^2$ by the diagonal linear automorphism with eigenvalues $e^{\frac{2\pi i}{n}}, \, e^{\frac{2\pi i q}{n}}$ with $(q,n)=1$. This singularity is called of \emph{type} $\frac{1}{n}(1,q)$, or $\frac{q}{n}$ for short.

\item  The exceptional divisor $E$ on the minimal resolution of a cyclic quotient singularity is given by a Hirzebruch-Jung string (see e.g., \cite[Chapter III, Section 5]{BHPV}).
A \emph{Hirzebruch-Jung string} (HJ-string, for short) is a union $\widetilde{E}:=\cup^k_i E_i$ of smooth
rational curves $E_i$ such that:
\begin{itemize}
\item $E^2_i = -b_i  \leq -2 \textrm{ for all } i$,
\item $E_iE_j=1$  if
$\mid i-j \mid =1$,
\item $E_iE_j=0$ if  $\mid i-j \mid \geq 2$,
\end{itemize}
where the $b_i$'s are given by the
\emph{continued fraction} associated to $\frac{1}{n}(1,q)$. Indeed, by the formula:
\[ \frac{n}{q}=b_1-\frac{1}{b_2-\frac{1}{...-\frac{1}{b_k}}}.
\]

\item A product-quotient surface comes together with two
isotrivial fibrations. Let us consider one of them: $\pi_2: X
\rightarrow C_2/G$. Take any point $b \in C_2/G$, and let $F$
denote the fibre of $\pi_2$ over $b$. Then (see \cite[Theorem 2.1]{S96}):
\begin{itemize}

\item The reduced structure of $F$ is the union of an irreducible smooth curve $Y$, called the \emph{central component} of $F$, and either none or at least two mutually disjoint HJ-strings, each one meeting $Y$ at one point. These strings are in one-to-one correspondence with the branch points of $C_1 \rightarrow (C_1/H)$, where $H \subset G$ is the stabilizer of $b$.
\item The intersection of a string with $Y$ is transversal and it takes place at only one of the end component of the string.

\end{itemize}

\item There are formulae for calculating the self intersection of the canonical divisor of a product-quotient surface.
\begin{equation*}\label{kiso}
 K^2_X\geq\frac{8(g(C_1)-1)(g(C_2)-1)}{| G |} +\sum_{x \in Sing (X)} h_x,
\end{equation*}
where  $h_x$ depends on the type of singularity at $x$ and the equality holds if $g(C_i)\geq 2$, for at least one value of $i\in\{1,2\}$.
If $x$ is a cyclic quotient singularity of type $\frac{1}{n}(1,q)$
then:
\[ h_x:= 2 - \frac{2+q+q'}{n}- \sum^k_{i=1}(b_i-2),
\]
where $q' \in \{1,\dots ,n-1\}$ is such that $qq' \equiv 1$ $mod \ n$, and
$[b_1, \dots ,b_k]$ is the continued fraction associated to $\frac{1}{n}(1,q)$. 

\item Finally, notice that $X$ is in general not a minimal model. Indeed, we will treat mostly examples with $X$ not minimal.

\end{enumerate}

\end{rem}


\subsection{Hodge Structure of $X$}\label{sec: Hodge structures}

Sometimes it will be useful to keep track of the action of an
automorphism on each curve separately; that is why in the
following we assume $G\simeq \ZZ/n\ZZ \simeq \langle x \rangle
\hookrightarrow \Aut(C_1)$ $x \mapsto {\bf g_1}$, and $G
\hookrightarrow  \Aut(C_2)$, $x \mapsto {\bf g_2}$. We have
$G=\langle {\bf g_1}\times {\bf g_2}\rangle\subset \Aut(C_1\times C_2)$, and we write $(C_1 \times
C_2)/({\bf g_1}\times {\bf g_2})$ for  $(C_1 \times
C_2)/(\langle{\bf g_1}\times {\bf g_2}\rangle)$.

We shall describe the Hodge structure of the product-quotient surface $X$ whose quotient model is $(C_1\times C_2)/G$.
The Hodge numbers of $X$ are determined by the action of $G$ on the cohomology of $C_1$ and $C_2$.

As in the previous section we denote by $\alpha_i$, $i=1,\ldots n$ the dimension of the
eigenspace $H^{1,0}(C_1)_{\zeta_n^i}$ with eigenvalue $\zeta^i$ w.r.t. the action of ${\bf g_1}$ and by $\beta_i$, $i=1,\ldots n$,
the dimension of the eigenspace $H^{1,0}(C_2)_{\zeta_n^i}$ with eigenvalue $\zeta_n^i$. w.r.t. the action of ${\bf g_2}$.

By \cite[Satz 1]{F71} we have

\[ H^0(X, \Omega^i_X) \simeq H^0(C_1\times C_2, \Omega^i_{C_1 \times C_2})^G.
\]
with $i=0,1,2$

Thus by the K\"unneth formula, we have

\begin{itemize}
\item $H^{0,0}(X)=H^{0,0}(C_1\times C_2)^G=H^{0,0}(C_1)\otimes H^{0,0}(C_2)$;

\item $H^{1,0}(X)=H^{1,0}(C_1\times C_2)^G=H^{1,0}(C_1)_{id}\otimes H^{0,0}(C_2)\oplus H^{0,0}(C_1)\otimes H^{1,0}(C_2)_{id}$,\\ in particular $h^{1,0}(X)=h^{0,1}(X)=\alpha_n+\beta_n$;
\item $H^{2,0}(X)=H^{2,0}(C_1\times C_2)^G=\sum_{i=1}^n H^{1,0}(C_1)_{\zeta_n^i}\otimes H^{1,0}(C_2)_{\zeta_n^{n-i}}$,\\ in particular $h^{2,0}(X)=h^{0,2}(X)=\sum_{i=1}^n\alpha_i\beta_{n-i}$.
\end{itemize}

In order to find $h^{1,1}(X)$, one has to know the number and the type of singularities of $(C_1\times C_2)/G$. Indeed, the desingularization of $(C_1\times C_2)/G$ introduces some exceptional divisors which increase the Picard number of the surface and thus $h^{1,1}(X)$.

Here we describe the structure of $H^{1,1}(X)$ starting from the description of the action of ${\bf g_i}$ on $C_i$, $i=1,2$.

Let us consider the set of points of $C_1$ (resp. $C_2$) with a non trivial stabilizer w.r.t. the action of ${\bf g_1}$ (resp. ${\bf g_2}$). Let us denote by $a_{i,h}$ (resp. $b_{i,h}$) the number of points on $C_1$ (resp. $C_2$) whose stabilizer has order $h$ and such that the local action of ${\bf g_1}^{|{\bf g_1}|/h}$ (resp. ${\bf g_2}^{|{\bf g_2}|/h}$) is $\zeta_{h}^i$.
Let $P\in C_1$ (resp. $Q\in C_2$) be a point with stabilizer of order $h$ (resp. $k$) where local action of ${\bf g_1}$ (resp. ${\bf g_2}$) is $\zeta_h^i$, $i \in \{1\ldots h\}$ (resp. $\zeta_k^j$, $j\in\{1\ldots k\})$. We assume $h$ (resp. $k$) are divisors of $|{\bf g_1}|$ different from $|{\bf g_1}|$.
We will say that $P\times Q$ is a point of type $((i,h),(j,k))$ and clearly the number of these points is $a_{i,h}b_{j,k}$. The stabilizer of these points has order $d_{(h,k)}:=\gcd(h,k)$ and the orbit of $P\times Q$ w.r.t. ${\bf g_1}\times{\bf g_2}$ contains $|{\bf g_1}|/d_{(h,k)}$ distinct points. In particular we observe that if $d_{(h,k)}=1$ the stabilizer of $P\times Q$ is empty and its orbit contains exactly $|{\bf g_1}|$ points and if $i=j=|{\bf g_1}|$, then the stabilizer of $P\times Q$ is $\ZZ/|{\bf g_1}|\ZZ$ and the orbit of $P\times Q$ consists only of $P\times Q$.
Let us now consider the quotient model $(C_1\times C_2)/G$ and the image $\overline{P\times Q}$ of the point $P\times Q$ for the quotient map: the point $\overline{P\times Q}$ is smooth if $d_{(h,k)}=1$; otherwise it is a singular points of type $\frac{1}{d_{(h,k)}}(1,q)$ where $q$ can be computed as in \cite[Prop. 5.3]{BHPV}. The images of the $a_{i,h}b_{j,k}$ points of type $((i,h),(j,k))$ under the quotient map consists of $a_{i,h}b_{j,k}d_{h,k}/|{\bf g_1}|$ points.\\

Now we assume that $n=p$ is a prime number: then $h$ (resp. $k$) is necessarily 1, hence $a_{i,h}=a_i$ (resp. $b_{j,k}=b_j$) with the notation of Section \ref{sec.curves}. The singular point $\overline{P\times Q}\in (C_1\times C_2)/G$ is of type $\frac{1}{p}(1,ij^{-1})$, where $ij^{-1}$ is computed in$\mod p$. Moreover, the orbit of every point with non trivial stabilizer consists only of one point and thus we have exactly $a_ib_j$ singular points on $(C_1\times C_2)/G$ which are of  type $\frac{1}{p}(1,ij^{-1})$.
\\ Once one knows the number and the type of singularities of $(C_1\times C_2)/G$ one can easily compute $h^{1,1}(X)$, recalling that every singularity we introduce a HJ-string, by Remark \ref{rem: PQ properties} (2),  and that
\begin{equation*}
 H^{1,1}(C_1\times C_2)^G=\left(H^{1,1}(C_1)\otimes H^{0,0}(C_2)\right)\oplus \left(\oplus_i \left(H^{1,0}(C_1)_{\zeta_n^i}\otimes H^{0,1}(C_2)_{\zeta_n^{n-i}}\right)\right)\oplus
 \end{equation*}
\begin{equation*}
\oplus
\left(\oplus_i \left(H^{0,1}(C_1)_{\zeta_n^i}\otimes H^{1,0}(C_2)_{\zeta_n^{n-i}}\right)\right)\oplus \left(H^{0,0}(C_1)\otimes H^{1,1}(C_2)\right),
\end{equation*}

In particular, since $H^{0,1}(C_j)_{\zeta_n^i}\simeq \overline{H^{1,0}(C_j)_{\zeta_n^{n-i}}}$, $j=1,2$,  $h^{1,1}(C_1\times C_2)^G=2(1+\sum_{i=1}^n\alpha_i\beta_i).$

The following proposition recap the results proved in this section.
\begin{prop}\label{prop: hodge numbers X} Let $C_i$, $i=1,2$ be a curve with an automorphism ${\bf g_i}$ of order $n$. Let $\alpha_l:=dim(H^{1,0}(C_1)_{\zeta_n^l})$ and $\beta_m:=dim(H^{1,0}(C_2)_{\zeta_n^m})$. Let $a_{(i,h)}$, $b_{(j,k)}$ and $d_{(h,k)}$ as above. Note that $i$ and $j$ are invertible in $\Z/d_{(h,k)}\Z$ and we denote by $r_i$ and $s_j$ their inverse. Let $X$ be the minimal resolution of the quotient surface $(C_1\times C_2)/({\bf g_1}\times {\bf g_2})$. Then the Hodge numbers of $X$ are:
$$h^{0,0}(X)=1, \ \ h^{1,0}(X)=\alpha_n+\beta_n,\ \ h^{2,0}(X)=\sum_{i=1}^n\alpha_i\beta_{n-1}$$
$$h^{1,1}(X)=2(1+\sum_{i=1}^n\alpha_i\beta_i)+\sum_{h,k}\sum_{i=1}^{n-1}\sum_{j=1}^{n-1}(a_{i,h}b_{j,k}d_{(h,k)}/n)k(i,j)$$ where:\\ $\bullet$ we pose $k(i,j)=0$ if $d_{(h,k)}=1$; \\
$\bullet$ $k(i,j)$ is the number of curves introduced by a singularity of type $\frac{1}{d_{(h,k)}}(1,q)$;\\
$\bullet$ $q:=i s_j\in\Z/d_{(h,k)}\Z$ (or equivalently $q:=jr_i\in\Z/d_{(h,k)}\Z$).
\end{prop}
\begin{rem}{\rm We recall that $\alpha_i$ and $\beta_i$ are uniquely determined by $a_i$ and $b_i$ and viceversa, thus the Hodge numbers of $X$ depend only on the set of values $\{\alpha_i\ \beta_i\}$, $i=1,\ldots, n$ (or equivalently on the set of values $\{a_i\ b_i\}$, $i=1,\ldots, n$).}\end{rem}

\subsection{The minimal model $S$}

As observed in the last point of Remark \ref{rem: PQ properties} the surface $X$ is in general non minimal. Let us denote by $S$ the minimal model of $X$. Since $h^{i,0}$ are birational invariant, the Hodge numbers $h^{0,0}$, $h^{1,0}$, $h^{2,0}$ of the product-quotient $X$ coincide with the ones of its minimal minimal model $S$. More complicate is the computation of $h^{1,1}(S)$.

In order to determine the minimal model $S$ of a product-quotient $X$, we have to find all the $(-1)$-curves on $X$. In the cases we will treat the $(-1)$-curves are central components of reducible fibres of one or both isotrivial fibrations of $X$. Then, after contractions of these, the $(-1)$-curves could be images of some divisors in the HJ-strings. Note that this is in general not true, see e.g., \cite{BP11}. A quick method to calculate the self intersection of the central components is given in \cite{P10}. We shall recall it.

\begin{defin}
We say that a reducible fibre $F_1$ of $\pi_2 \colon X \lr
C_2/G$
 is of type $\big( \frac{q_1}{n_1}, \ldots, \frac{q_r}{n_r} \big)$ if
it contains exactly $r$ $HJ$-strings $\widetilde{E}_1, \ldots, \widetilde{E}_r$, where
each $\widetilde{E}_i$ is of type $\frac{1}{n_i}(1, q_i)$. The same definition
holds for a reducible fibre $F_2$ of $\pi_1 \colon S \lr
C_1/G$.
\end{defin}
\begin{prop}\cite[Proposition 2.8]{P10} \label{prop: selfinter central}
Let $F_1$ be of type $\big( \frac{q_1}{n_1},  \ldots,
\frac{q_r}{n_r} \big)$ and let $Y_1$ be its central component.
Then
\begin{equation*} \label{selfint-1}
(Y_1)^2 = - \sum_{i=1}^r \frac{q_i}{n_i}.
\end{equation*}
If $F_2$ is of type $(\frac{q'_1}{n_1},  \ldots,
 \frac{q'_r}{n_r})$ then
$(Y_2)^2 = - \sum_{i=1}^r \frac{q_i'}{n_i}.$
\end{prop}

In the following example we construct two surfaces, $S_1$ and $S_2$, both of them are the minimal model of quotients of $D_3\times D_3$ by a diagonal action of $\Z/3\Z$ but the action of this group is different in the two cases. As a consequence the minimal resolution of one quotient has an infinite number of $(-1)$-curves, the minimal resolution of the other has no $(-1)$-curves and the minimal models of these two resolutions are totally different: one of them is a rational surface, one is a K3 surface.

\begin{example}\label{example: quotient hodge numbers 0,i} Let us consider the product of two elliptic curves $D_3\times D_3$ and its automorphisms $\delta_3\times \delta_3^1$ and $\delta_3\times \delta_3^2$. We will denote by $X_i$ the minimal resolution of $(D_3\times D_3)/(\delta_3\times \delta_3^i)$ and by $S_i$ its minimal model. We recall that $\delta_3^i$ acts on $H^{1,0}(D_3)$ as the multiplication by $\zeta_3^i$. We obtain $h^{0,0}(X_i)=h^{0,0}(S_i)=1$, $h^{1,0}(X_i)=h^{1,0}(S_i)=0$ for $i=1,2$; $h^{2,0}(X_1)=h^{2,0}(S_1)=0$ and $h^{2,0}(X_2)=h^{2,0}(S_2)=1$.\\ Now we compute $h^{1,1}(X_i)$: on $(D_3\times D_3)/(\delta_3\times \delta_3)$ there are 9 singularities of type $\frac{1}{3}(1,1)$ and on $(D_3\times D_3)/(\delta_3\times \delta_3^2)$ there are 9 singularities of type $\frac{1}{3}(1,2)$. The resolution of a point of type $\frac{1}{n}(1,1)$ introduces 1 curve of self intersection $-n$, and thus $h^{1,1}(X_1)=9+dim\left(H^{1,1}(D_3\times D_3)^{\delta_3\times \delta_3}\right)=13$. The desingularization of singularities of type $\frac{1}{n}(1,n-1)$ introduces $n-1$ rational curves with self intersection $-2$, whose dual diagram is $A_{n-1}$, and thus $h^{1,1}(X_2)=2\cdot 9+dim\left(h^{1,1}(D_3\times D_3)^{\delta_3\times \delta_3^2}\right)=20$.\\
With the same method of the Section \ref{sec: equation} we obtain that an equation for $(D_3\times
D_3)/(\delta_3\times \delta_3^i)$ is given by
$y^2=x^3+(v^2-1)^{2i}$, which is the equation of an elliptic
fibration over $\mathbb{P}^1_{[v]}$. If $i=1$, this elliptic
fibration is a rational elliptic fibration (this depend on the
degree of $(v^2-1)^{2i})$. By the standard theory of elliptic
fibration (cf. \cite{M}), it has 3 reducible fibers (over $1$, $-1$
and $\infty$)
of type $IV$ (coming from the contraction of the central components of the reducible fibers on $X_1$). The rank of the Mordell--Weil group of this elliptic fibration is $10-2-6=2$. This implies that there are infinite sections of this elliptic fibration, and thus an infinite number of $(-1)$-curves. Since the minimal elliptic fibration is obtained by contractions of $X_1$, we obtain that $X_1$ contains an infinite number of $(-1)$-curve. A minimal model, $S_1$, of $X_1$ is birational to $\mathbb{P}^2$, so $h^{1,1}(S_1)=1$.\\
The elliptic fibration $y^2=x^3+(v_2^2-1)^{4}$ has a K3 surface as minimal model. The reducible fibers of this elliptic fibration are 3 fibers of type $IV^*$ each of them consists of the central component (which is a $(-2)$-curve, by Proposition \ref{prop: selfinter central}) and of three copy of $A_2$. So $X_2$ coincides with $S_2$ and $h^{1,1}(S_2)=20$. This K3 surface was constructed in \cite{SI}.
\end{example}

\subsection{Automorphisms of $S$ and quotient surfaces}\label{subsect: autom S}
By construction the surface $(C_1\times C_2)/({\bf g_1}\times {\bf g_2})$ always admits an automorphism of order $n$ induced by $\id \times {\bf g_2}\in \Aut (C_1\times C_2)$ (or equivalently by ${\bf g_1}\times id\in \Aut(C_1\times C_2)$). This automorphism lifts to an automorphism of $X$ and of $S$. Thus one can consider the quotient $\left(\left(C_1\times C_2\right)/\left( {\bf g_1}\times {\bf g_2}\right)\right)/\left(\id \times {\bf g_2}\right)$. Since $\langle {\bf g_1}\times {\bf g_2},\id\times {\bf g_2}\rangle=\langle {\bf g_1}\times \id, \id\times {\bf g_2}\rangle$, we have the following commutative diagram:
$$
\begin{array}{rccccc}
&&C_1\times C_2&&\\
&\swarrow_{n:1}&&\searrow^{n:1}&\\
\left(C_1\times C_2\right)/{\bf g_1}\times {\bf g_2}&&&&\left(C_1\times C_2\right)/\id\times {\bf g_2}\\
&\searrow^{n:1}&&\swarrow_{n:1}&\\
&&C_1/{\bf g_1}\times C_2/{\bf g_2}.&&
\end{array}
$$
This diagram lifts to the minimal resolution of all the surfaces we are considering, and so one obtains that the surface $X$ has a generically $n^2:1$ map to $C_1/{\bf g_1}\times C_2/{\bf g_2}$, and in particular to $\mathbb{P}^1\times \mathbb{P}^1$, if we assume $C_i/{\bf g_i}\simeq \mathbb{P}^1$. Moreover, the map $X\ra \mathbb{P}^1\times \mathbb{P}^1$ induces a rational $n^2:1$ map $S\dashrightarrow \mathbb{P}^1\times \mathbb{P}^1$.\\

We can explicitly describe the action of $\id \times {\bf g_2}$ on the cohomology of $X$ (we keep the assumption $C_i/{\bf g_i}\simeq \mathbb{P}^1$):
$\id \times {\bf g_2}$ acts as the identity on the spaces $H^{0,0}(X)$ and $H^{1,0}(X)$ . The invariant subspace of $H^{2,0}(X)$ under $\id \times {\bf g_2}$ is the image of the space $H^{1,0}(C_1)^{\bf g_1}\otimes H^{1,0}(C_2)^{\bf g_2}$ and thus has dimension $\alpha_n\beta_n=0$. As we saw in Section \ref{sec: Hodge structures}, the space $H^{1,1}(X)$ splits into two parts: the image of $H^{1,1}(C_1\times C_2)^{{\bf g_1}\times{\bf  g_2}}$ and a direct summand, say $R$, which comes from the resolution of the singularities of $(C_1\times C_2)/({\bf g_1}\times {\bf g_2})$. Hence, $H^{1,1}(X)^{\id \times {\bf g_2}}$ splits into the direct sum of the image of $(H^{0,0}(C_1)\otimes H^{1,1}(C_2))\oplus (H^{1,0}(C_1)^{\bf{g_1}}\otimes H^{0,1}(C_2)^{{\bf g_2}})\oplus (H^{0,1}(C_1)^{{\bf g_1}}\otimes H^{1,0}(C_2)^{{\bf g_2}})\oplus (H^{1,1}(C_1)\otimes H^{0,0}(C_2))$ and $R^{\id \times\bf g_2 }$. The dimension of the first term is $2+2\alpha_n\beta_n=2$. We note that, if every point in the branch locus of $C_i$, $i=1,2$ is of total ramification, then the action of $\id\times {\bf g_2}$ is the identity on $R$, and one finds
$H^{1,1}(X)^{\id \times {\bf g_2}}=(H^{0,0}(C_1)\otimes H^{1,1}(C_2))\oplus (H^{1,1}(C_1)\otimes H^{0,0}(C_2))\oplus (R\otimes \C)$.

\section{K3 surfaces}\label{Sec K3 sur}
This section is devoted to the construction of K3 surfaces $S$, which are minimal model of product-quotient surfaces.

We recall that by definition a \emph{K3 surface} $S$ has $h^{1,0}(S)=0$ and trivial canonical bundle. The Hodge numbers of $S$ are uniquely determined by these properties and are $h^{0,0}(S)=h^{2,0}(S)=1$, $h^{1,0}(S)=0$, $h^{1,1}(S)=20$.


\subsection{Product-quotient surfaces with $p_g=1$ and $q=0$} \label{subset: quotient Zp with pg1 q0}

Let $S$ be the minimal model of a product-quotient $X$ with quotient model $(C_1\times C_2)/({\bf g_1} \times {\bf g_2})$. If $S$ is a K3 surface, then $q(S)=h^{1,0}(S)=0$ and $p_g(S)=h^{2,0}(S)=1$. Since $h^{i,0}$ are birational invariants, $h^{1,0}(X)=0$, $h^{2,0}(X)=1$. Therefore, by Proposition \ref{prop: hodge numbers X}, $\alpha_n+\beta_n=0$ and $\sum_{i=1}^n\alpha_i\beta_{n-i}=1$. In particular
\begin{equation}\label{eq: condition pg1 q0}
\alpha_n=\beta_n=0\mbox{ and there exists }i\in\{ 1\ldots n\}\mbox{ such that }\alpha_i=\beta_{n-i}=1,\ \ \alpha_j\beta_{n-j}=0\mbox{ if }j\neq i.\\
\end{equation}
Condition \eqref{eq: condition pg1 q0} is divided in two: a condition on each factor of the product, namely the action of ${\bf g_i}$ on $H^{1,0}(C_i)$; and a condition on the whole product, namely the action of ${\bf g_1} \times {\bf g_2}$ on $H^{2,0}(C_1\times C_2)$. \\

Now we shall give a procedure to construct product-quotient surfaces  with $p_g=1$ and $q=0$. Since, in general, this algorithm requires excessively long calculations, many of them are performed using the \verb|MAGMA| script in the Appendix. While describing the procedure we also indicate which part of the program does what.

First we search for pairs $(C_1,{\bf g_1})$ which satisfy the following two conditions: $\alpha_n=0$, or equivalently $C_i/{\bf g_i}\simeq \mathbb{P}^1$; and there exists an index $0<j<n$ such that $\alpha_j=1$. The first condition is obtained by giving spherical systems of generators for $G$ (see Theorem \ref{theo: Riemann}). The second one is obtained exploiting the Chevalley--Weil formula and calculating the rotation angles as in Proposition \ref{prop: rotation consts}. Analogous conditions must hold for a second pair $(C_2,{\bf g_2})$.

Second we couple the curves by requiring that there exists an $i\in \{ 1\ldots n\}$ such that $\alpha_i=\beta_{n-i}=1$ and $\alpha_j\beta_{n-j}=0$ if $j\neq i$. \\

Let us assume that $n=p$ a prime. The curves $C_i$, admitting an automorphism ${\bf g_i}$ of order $p$ such that $C_i/{\bf g_i}\simeq \mathbb{P}^1$ and $\alpha_j=1$ for some $j\in \{ 1,\ldots, p-1\}$, have genus at most $(p-1)^2$ by Corollary \ref{cor: max genus}. The number of the ramification points of the cover $C_i\ra C_i/{\bf g_i}$ is at most $2p$. This implies that the number of curves with these properties is finite. These curves are classified by the  \verb|MAGMA| program given in the Appendix. The calculation can be found in the function \verb|Surfacesp|. The function calculates all the partitions of  all the numbers from 3 up to $2p$, giving all the admissible ramification data $(a_1,\ldots , a_{n-1})$ of the coverings $C_i\ra C_i/g_i\simeq \mathbb{P}^1$. Afterwards it evaluates the $\alpha_j$, as in Proposition \ref{prop: alpha ai}, and lists only the ones with $\alpha_j=1$ for at least one $j\in \{ 1,\ldots, p-1\}$.

\begin{rem}\label{rem: existence of curves p}{\rm We observe that for every prime $p$ there exists at least one curve with the required properties, the one with $a_{p-1}=2$ and $a_{(p+1)/2}=1$ (cf. Example \ref{ex: curve Dp}).}\end{rem}

The condition $\alpha_j\beta_{n-j}=0$ if $j\neq i$ implies that the list $(\alpha_1,\ldots, \alpha_{p-1},\beta_1,\ldots \beta_{p-1})$ contains at least $p-2$ zeros. This condition is verified in the function \verb|MaybeSur1| in the Appendix. Then, the function \verb|TheSur| tests if a surface given by \verb|MaybeSur1| has also $p_g=1$.

\begin{prop}\label{prop: quotient Zp pg1 q0}
There exists a finite number of surfaces $S$ which are the minimal model of $(C_1 \times C_2)/({\bf g_1} \times {\bf g_2})$ with $|({\bf g_1}\times {\bf g_2})|=p$ and $p_g(S)=1$, $q(S)=0$.
\end{prop}
\proof For a given $p$, the number of surfaces $S$ is finite, since the numbers of pairs $(C_i,{\bf g_i})$ is so (cf. Corollary \ref{cor: max genus}).

The automorphism $\id\times {\bf g_2}$ induces an automorphism on $S$ which acts non trivially on $H^{2,0}(S)$, see Section \ref{subsect: autom S}. In order to give a bound for $p$ we prove that there exists no a surfaces $Z$ with $p_g(Z)=1$, $q(Z)=0$ and an automorphism of order $p>19$ acting non trivially on $H^{2,0}(Z)$. \\
Let $Z$ be a minimal surface with $h^{2,0}(Z)=1$ which admits an automorphism $\sigma$ of order $p$ and let $H^{2}(Z)_{\zeta_p^i}$ be the eigenspace of the eigenvalue $\zeta_p^i$ for the action of $\sigma$, $i=0,\ldots, p-1$. The dimension of $H^{2}(Z)_{\zeta_p^i}$ does not depend on $i$ if $i\neq 0$. Thus if there exists $i\neq 0$, $i\in\{ 1 \ldots p-1 \}$ such that $dim(H^{2}(Z)_{\zeta_p^i})\geq 1$, $b_2(Z)\geq p-1$. \\
Since $h^{2,0}(Z)=1$, $K_Z^2\geq 0$. By $h^{1,0}(Z)=0$ follows that $\chi(Z)=2$ and that $e(Z)\leq 24$ by Noether equality. So $b_2(Z)\leq 22$.\\
Since there exists no a surface $Z$ with $p_g(Z)=1$ $q(Z)=0$ and $b_2(Z)>22$, there exists no a surfaces $Z$ with $p_g(Z)=1$, $q(Z)=0$ and an automorphism of order $p>19$ acting non trivially on $H^{2,0}(Z)$.\endproof

\begin{rem}\label{rem: p leq19}{\rm Since a minimal surface with $p_g=1$, $q=0$ can not admit an automorphism of order $p>19$, there exists no pairs $(C_1\times C_2, {\bf g_1}\times {\bf g_2})$ such that $|{\bf g_1}\times {\bf g_2}|=p>19$, $dim (H^1(C_1\times C_2)^{{\bf g_1\times g_2}})=0$ and $dim (H^{2,0}(C_1\times C_2)^{{\bf g_1\times g_2}})=1$.}\end{rem}

\begin{lem}\label{lemma: h0(C)} Let $X$ be the minimal resolution of $(C_1\times C_2)/G$. Let $F_i$ be the fiber of $\pi_j:X\ra C_j/G$, $\{i,j\}=\{1,2\}$. We recall $F_i\simeq C_i$. If $q(X)=0$, the linear systems $|F_i|$ on $X$, $i=1,2$ are complete and of dimension 1.\end{lem}
\proof Since $F_i^2=0$, we have the short exact sequence 
$$0\ra \mathcal{O}(K_X-F_i)\ra \mathcal{O}(K_X)\ra \mathcal{O}_{F_i}(K_{F_i})\ra 0.$$
This induces the long exact sequence in cohomology
$$H^1(X,\mathcal{O}(K_X))\ra H^1(F_i,\mathcal{O}(K_{F_i}))\ra H^2(X,\mathcal{O}(K_X-F_i))\ra H^2(X,\mathcal{O}(K_X))\ra 0.$$
By Serre duality $\dim (H^1(X,\mathcal{O}(K_X)))=q(X)=0$ and  $H^2(X,\mathcal{O}(K_X-F_i))=H^0(X,\mathcal{O}(F_i))$. Therefore  $\dim(H^0(X,\mathcal{O}(F_i)))=\dim(H^1(F_i,\mathcal{O}(K_{F_i})))+\dim(H^2(X,\mathcal{O}(K_X)))=\dim(H^0(F_i,\mathcal{O}))+1=2.$
\endproof

\begin{prop}\label{prop: moduli=r1+r2-6}  Let $X$ be the minimal resolution of $(C_1\times C_2)/G$. If $q(X)=0$, the dimension of the family of product-quotient surfaces of $(C_1\times C_2)/G$ is $r_1+r_2-6$, where $C_i\ra C_i/G$ is ramified in $r_i$ points.
\end{prop}
\proof
By Teichm\"uller theory the number $\eta$ of parameters of the family of product-quotient surfaces is less then or equal to $r_1+r_2-6$. For simplicity we assume $r_2=3$. If $\eta<r_1-3$, then there exists a positive dimensional family of curves isomorphic to $C_1$ which induces a family of isotrivial fibrations on $X$, whose fibers are isomorphic to $C_1$. Since $q(X)=0$, $Pic(X)=NS(X)/{\rm Tors}\simeq \Z^N/{\rm Tors}$ for a certain positive integer $N$. By Lemma \ref{lemma: h0(C)}, the linear system $|F_1|$ is complete and of dimension 1. Since $Pic(X)$ is discrete there is no positive dimensional family of such linear systems.  Therefore $\eta=r_1-3$.
\endproof
There is a finite number of surfaces as described in Proposition \ref{prop: quotient Zp pg1 q0}. These surfaces are given by the program \verb|Surfacesp|. However
the number of permutations of the ramification points increases
too rapidly with the growth of $p$ for a computation in a short
time. Since our aim was the construction of K3 surfaces, and we
know the dimension of the families we are searching for, we
wrote another program with a fixed number of ramification
points, and hence with fixed dimension of the family. The program \verb|t1t2PtsSurfaces| in the Appendix --
given a cyclic group $G$ of order $p$ or $2p$ and the numbers
$ti$, $i=1,2$, of ramification points of $C_i \rightarrow C_1/G
\simeq \PP^1$ -- returns a list of product-quotient surfaces $X$
with $p_g(X)=1$ and $q(X)=0$, as well as the singularities of $X$.

\begin{rem}\label{rem: complete list} If $|{\bf g_1}|=2p$, then by Remark \ref{rem: max genus 2p} there exists a finite list of curves with at least one $\alpha_j=1$ and $\zeta_{2p}^j$ is a primitive $2p$-root of unity. If the action of ${\bf g_1}$ on $C_1$ is of this type, the same must be true for the action of ${\bf g_2}$ on $C_2$, thus we have a finite list for $(C_1\times C_2)/({\bf g_1} \times {\bf g_2})$. Hence, we obtain a complete classification of such surfaces as in Proposition \ref{prop: quotient Zp pg1 q0}. Otherwise, if we assume that the action of ${\bf g_1}$ on $C_1$ is such that $\alpha_j=1$ and $\zeta_{2p}^j$ is a either a primitive $p$-root of unity or $(-1)$, then the same must be true for the action of ${\bf g_2}$ on $C_2$. In this case we can not construct a complete list of the curves $C_1$ and $C_2$, since we have no an upper bound for their genera, and so for the number of ramification points of the map $f_i:C_i\ra C_i/{\bf g_i}\simeq \mathbb{P}^1$. Anyway, if we fix the maximal number $n$ of ramification points for $f_i$, then we obtain a finite list of curves $C_i$ and thus a finite number of  surfaces $(C_1\times C_2)/({\bf g_1}\times {\bf g_2})$ as in Proposition \ref{prop: quotient Zp pg1 q0}. Our aim is to construct K3 surfaces, so the bound on $n$ depends on the dimension of the moduli space of K3 surfaces. More precisely, the moduli space of projective K3 surfaces has dimension 19, so the sum of the ramification points of $f_1$ and $f_2$ can not be grater then 25.
\end{rem}


\subsection{K3 surfaces}
Let $S$ be the minimal model of $(C_1\times C_2)/({\bf g_1} \times {\bf g_2})$.
Let us assume $\langle {\bf g_1} \times {\bf g_2} \rangle\simeq
\Z/n\Z$. If $S$ is a K3 surface, then by definition $p_g(S)=1$ and
$q(S)=0$. Therefore the K3 surfaces obtained as minimal model of
$(C_1\times C_2)/({\bf g_1}\times {\bf g_2})$ are among the ones
listed in Section \ref{subset: quotient Zp with pg1 q0}.
In order to prove that $S$ is a K3 surface one has to verify that the canonical bundle is trivial.
\begin{lem}\label{lem: trivial canonical bundle}
Let $Z$ be a surface obtained contracting $-K_X^2$ $(-1)$-curves on $X$. We recall that $p_g(Z)=1$ and $q(Z)=0$. Let $F_1$ be the class of the fiber of the fibration $\pi_2:X \ra C_2/\langle {\bf g_2} \rangle$ and
$F_2$ be the class of the fiber of the fibration $\pi_1:X\ra
C_1/\langle{\bf  g_1} \rangle$. Let $E$ be the sum  of all the
exceptional divisors of the blow up $X\ra Z$. If $(K_X-E)F_1=0$
and $(K_X-E)F_2=0$, then $K_Z$ is trivial. In this case $Z$ is minimal.
\end{lem}
\proof
We shall denote by $P_k$ the singular points of $(C_1\times C_2)/({\bf g_1}\times {\bf g_2})$, and by $A_{j,k}$ the exceptional curves of the blow up $\pi:X\ra(C_1\times C_2)/({\bf g_1} \times {\bf g_2})$ over $P_k$.\\
Let $D$ be an effective divisor on $X$, then $D=\lambda_1F_1+\lambda_2F_2+\lambda_3B+\sum_{j,k} \lambda_{j,k} A_{j,k}$ with $\lambda_i, \lambda_{j,k}\geq 0$ and $BF_1>0$, $BF_2>0$. If $DF_1=0$ and $DF_2=0$, then $\lambda_1=\lambda_2=\lambda_3=0$.
For every $k$, $A_{j,k}$  is a HJ-string, hence $DA_{j,k}=0$ for every $j$ and $k$. Therefore $DF_1=DF_2=0$ give a homogeneuos linear system in $\lambda_{j,k}$. 
The corresponding matrix is a diagonal block matrix and each block is invertible, being associated to the resolution of the quotient singularity $P_k$. Thus $D=0$.\\
The divisor $K_X-E$ is an effective divisor being the pullback of the canonical divisor of $Z$, which has $p_g(Z)=1$. Applying the previous result to $D=K_X-E$ we obtain $K_X-E=0$ hence $K_Z$ is trivial.\endproof

Necessary conditions to obtain a K3 surface $S$ as minimal model of the minimal resolution $X$ of the quotient $(C_1\times C_2)/({\bf g_1}\times {\bf g_2})$ are the following:
\begin{enumerate}
\item $h^{i,0}(X)=1$, for $i=0,2$; \item $h^{1,0}(X)=0$; \item
there are exactly $-K_X^2$ $(-1)$-curves on $X$.
\end{enumerate}
Thus, in order to classify the K3 surfaces $S$, we list the surfaces $X$ satisfying the conditions (1), (2). If $|{\bf g_1}\times {\bf g_2}|=p$ is  prime number, as already said, this is done by the program \verb|Surfacesp|. \\

Next we consider the $(-1)$-curves which are either central components $Y$ of reducible fibers (we calculate $Y^2$ using Proposition \ref{prop: selfinter central}) or appear as image of curves in HJ-string after some contractions. In this way we find $-K_X^2$ $(-1)$-curves. After the contraction of all these curves we always obtain a surface which satisfy the condition of Lemma \ref{lem: trivial canonical bundle} and so a surface $Z$ with trivial canonical bundle. This implies that there are no other $(-1)$-curves on $Z$, which is thus the minimal model $S$ of $X$ and in particular $S$ is a K3 surface.

\subsection{Non-symplectic automorphisms}\label{sec: non sympl autom}
We saw in Section \ref{subsect: autom S} that every surface $S$ which is the minimal model of a product-quotient with ${\bf g_1}\times {\bf g_2}$ admits an automorphism induced by $\id \times {\bf g_2}$ which acts non trivially on $H^{2,0}(S)$. This means that if $S$ is a K3 surface, the automorphism induced on $S$ by $\id\times {\bf g_2}$ is a purely non-symplectic automorphism on $S$. Thus the surface $S$ admits a non-symplectic automorphisms of prime order.
\begin{defin} Let $W$ be a K3 surface. Let $\omega_W$ be a generator of $H^{2,0}(W)$. An automorphism ${\bf g}\in Aut(W)$ of order $n$ is called \emph{symplectic} if ${\bf g}(\omega_W)=\omega_W$, \emph{purely non-symplectic} if ${\bf g}(\omega_W)=\zeta_n^i\omega_W$ and $\zeta_n^i$ is a primitive $n$-root of unity.\end{defin}
We observe that an automorphism of prime order $p$ which is non-symplectic is purely non-symplectic.This type of automorphism are classified \cite{AST}. In this section we summarize the main results on non-symplectic automorphisms of prime order, which will be considered in the following.\\

For every prime number $2\leq p\leq 19$ there exists a K3 surface $W$ admitting a non-symplectic automorphism ${\bf g}$ of order $p$.
Let us assume $3\leq p\leq 19$. The fixed locus $Fix_{\bf g}(W)=\{w\in W \mbox{ such that }{\bf g}(w)=w\}$ consists of the disjoint union of $n$ isolated points and $k+1$ curves. At most one of the fixed curves has a positive genus, and we denote by $g(C)$ the genus of the curve with highest genus. Hence, the fixed locus consists of $n$ isolated points, $k$ rational curves and another curve $C$ with a possibly positive genus. For each prime number $3\leq p\leq 19$ there exists a finite number of possibilities for the fixed locus $Fix_{\bf g}(W)$, and the fixed locus is uniquely determined by the three invariants $(n,g( C ), k+1)$. The admissible choices for $(n,g( C ), k+1)$ are listed in \cite{AST}, where it is also proved that there exists a K3 surface with a non-symplectic automorphism of order $p$ with fixed locus determined by $(n,g( C ), k+1)$ for every admissible choice of the invariants.\\
More precisely, the invariants $(p,n,g(C),k+1)$ determine uniquely the two lattices $S_{(n,g( C ), k+1)}^p:=H^2(W,\Z)^g$ and $T_{(n,g( C ), k+1)}^p:=(H^2(W,\Z)^g)^{\perp}$ and a K3 surface admits a non-symplectic automorphism of order $p$ with fixed locus determined by $(n,g( C ), k+1)$ only if $S_{(n,g( C ), k+1)}^p$ is primitively embedded in its N\'eron-Severi group. This allows one to describe the family of K3 surfaces with a non-symplectic automorphism of order $p$ and a prescribed fixed locus in terms of the period map of K3 surfaces. We will denote by $\mathcal{M}^p_{(n,g( C ), k+1)}$ the family of K3 surfaces admitting a non-symplectic automorphism of order $p$ with fixed locus determined by $(n,g( C ), k+1)$. It has one connected component of dimension $\left(\rk (T_{(n,g( C ), k+1)}^p)/(p-1)\right)-1$. To give a more precise description of the moduli space of the K3 surfaces that admit a non-symplectic automorphism of order $p$ and a prescribed fixed locus, we consider the action of ${\bf g}$ on $T_{(n,g( C ), k+1)}^p\otimes\C$, which does not depend on the K3 surface considered. By definition ${\bf g}$ has no eigenvalue 1 on $T_{(n,g( C ), k+1)}^p\otimes\C$ and the decomposition in eigenspaces consists of $p-1$ equidimensional eigenspaces (of the eigenvalues $\zeta_p^i$, $i=1,\ldots p-1$). Let $(T_{(n,g( C ), k+1)}^p\otimes\C)_{\zeta_p}$ be the unique eigenspace such that $(T_{(n,g( C ), k+1)}^p\otimes\C)_{\zeta_p}^{2,0}\neq 0$. Set $\mathcal{B}:=\{z\in \mathbb{P}((T_{(n,g( C ), k+1)}^p\otimes\C)_{\zeta_p})\mbox{ such that }(z,z)=0,\ \ (z,\overline{z})>0\}$. The space $\mathcal{B}$ is a ball of dimension $(\rk (T_{(n,g( C ), k+1)}^p)/(p-1))-1$. Let $\Gamma:=\{\gamma\in {\rm O}(T_{(n,g( C ), k+1)}^p\otimes\C)\mbox{ such that }\gamma {\bf g}={\bf g}\gamma\}.$ The generic point of $\mathcal{B}/\Gamma$ corresponds to a K3 surface admitting a non-symplectic automorphism as required and there is a birational map between the space of such K3 surfaces and $\mathcal{B}/\Gamma$ (see \cite[Section 11]{DK}).\\

For a fixed prime number $3\leq p\leq 19$, there are some inclusions among the families $\mathcal{M}^p_{(n,g( C ), k+1)}$, for example $\mathcal{M}^3_{(9,0,6)}\subset \mathcal{M}^3_{(8,0,5)}$. All these inclusions are described in \cite{AST} and the maximal components are determined in \cite[Theorem 9.5]{AST}: \\
if $p=3$, there are three maximal components: $\mathcal{M}^3_{(3,-,0)}$, $\mathcal{M}^3_{(0,4,1)}$, $\mathcal{M}^3_{(0,5,2)}$ and for every admissible data $(n,g(C),k+1)$,  we have $\mathcal{M}^3_{(n,g(C),1)}\subset \mathcal{M}^3_{(0,4,1)}$ and $\mathcal{M}^3_{(n,g(C),k+1)}\subset \mathcal{M}^3_{(0,5,2)}$ if $k+1\geq 2$. The dimension of each family is $m=10-n$;\\
if $p=5,7,11$, there are two maximal components: $\mathcal{M}^p_{(n,-,0)}$,  and $\mathcal{M}^3_{(n,g(C),1)}$ and for every admissible data $(n,g(C),k+1)$,  we have $\mathcal{M}^3_{(n,g(C),k+1)}\subset \mathcal{M}^3_{(n,g(C),1)}$ if $k+1\neq 0$. The dimension of each family is
$m=(13-n)/(p-2)$ if $p=5,7$ and $m=(11-n)/(p-2)$ if $p=11$;\\
if $p=13,17,19$, there is only one (rigid) K3 surface admitting a non-symplectic automorphism of order $p$. Thus there is one maximal component (in fact one component of dimension 0) which is $\mathcal{M}^{13}_{(9,0,1)}$, $\mathcal{M}^{17}_{(7,-,0)}$, $\mathcal{M}^{19}_{(5,-,0)}$ respectively.

The K3 surfaces we are constructing are members of the families of K3 surfaces admitting a non-symplectic automorphism of prime order. In certain cases it turns out they are the general member of some of these families. In the following sections we will construct K3 surfaces $S$ and we will determine the fixed locus of the non-symplectic automorphism induced by $\id \times {\bf g_2}$ (or by some of its powers), and thus we identify on which component of the family of K3 surfaces with a non-symplectic automorphisms they lie.

The following remark is used to determine the fixed locus.
\begin{rem}\label{rem: induced automorphism} Let us assume $|{\bf g_2}|=p$ is a prime number. The automorphism of $S$ induced by
$\id\times {\bf g_2}$ fixes all the central components of all the reducible fibers of both the fibrations $\pi_1$ and $\pi_2$ and all the singular points of the HJ-strings introduced resolving the singularities of $(C_1\times C_2)/({\bf g_1}\times {\bf g_2})$. It is possible that the automorphism fixes some disjoint components of certain HJ-strings which do not meet any other fixed curves. It induces an automorphism ${\bf g}_S$ on $S$ which is non-symplectic of order $p$.
\end{rem}


\section{K3 surfaces which are minimal models of $(C_1\times C_2)/(\Z/p\Z)$}\label{sec: PQ with Zp}

The aim of this section is to list and to describe the K3 surfaces obtained as minimal model of $(C_1\times C_2)/({\bf g_1}\times {\bf g_2})$ with $|{\bf g_1}\times {\bf g_2}|=p$.
\begin{theo}\label{theo: PQ with Zp}
The K3 surfaces $S$ admitting a non-symplectic automorphism of odd prime order $p$ with fixed locus listed in the column $(n,g,k+1)$ of the \verb|Table 1| are all minimal models of product-quotient surfaces with group $\Z/p\Z$. Moreover, for each such surface $S$ we can choose $(C_2,{\bf g_2})\simeq (D_p,\delta_p)$ and the non-symplectic automorphism on $S$ is always induced by $\id\times\delta_p$.

Viceversa all the K3 surfaces which are minimal models of $(C_1\times C_2)/({\bf g_1}\times {\bf g_2})$, $|{\bf g_1}\times {\bf g_2}|=p$ admit a non-symplectic automorphism of order $p$ whose fixed locus is one of those listed in \verb|Table 1|.\end{theo}
\proof
Let us fix $p$. By Section \ref{sec: non sympl autom} we know the dimension $M_p$ of the maximal components of the family of K3 surfaces with a non-symplectic automorphim of order $p$. Since every K3 surface minimal model of a product-quotient with group $\Z/p\Z$ admits a non-symplectic automorphism of order $p$ (see Section \ref{subsect: autom S}), we can bound the number of moduli of the pairs $(C_1,{\bf g_1})$, $(C_2,{\bf g_2})$ by $M_p$, see also Proposition \ref{prop: moduli=r1+r2-6}.

The first step consists in listing product-quotient surfaces with $p_g=1$, $q=0$. This is done using the program \verb|t1t2PtsSurfaces| giving the group $G\simeq \Z/p\Z$ and the numbers $t1$, $t2$ such that $t1\geq 3$, $t2\geq 3$ and $t1+t2=m-6$, where $m\leq M_p$. Indeed, recall that $t1$ and $t2$ are the numbers of branching points of the projections $C_i\ra C_i/{\bf g_i}$ respectively. Then $t1-3$ and $t2-3$ are the moduli of the pairs $(C_i,{\bf g_i})$, see e.g. \cite{cat00}, and $m$ is the dimension of family of the product-quotient surfaces, by Proposition \ref{prop: moduli=r1+r2-6}.

Second step: For every product-quotient surface in the list one has to check if the minimal model is a K3 surface and has to calculate the fixed locus, determining $(n,g,k+1)$, of the induced automorphism. This is done for every entry in the list exactly as in Example \ref{subsec: example p=5}.

Every member of the family, $\mathcal{F}_{PQ}(C_1,C_2)$, of K3 surfaces minimal models of $(C_1\times C_2)/({\bf g_1}\times {\bf g_2})$, is also a member of a family $\mathcal{M}^p_{(n,g,k+1)}$. This implies that $\mathcal{F}_{PQ}(C_1,C_2)\subset \mathcal{M}^p_{(n,g,k+1)}$. For every $(n,g,k+1)$, there exists always a choice of $(C_1,{\bf g_1})$ and $(C_2,{\bf g_2})$ such that $\dim(\mathcal{F}_{PQ}(C_1,C_2))= \dim(\mathcal{M}^p_{(n,g,k+1)})$ which implies $\mathcal{F}_{PQ}(C_1,C_2)= \mathcal{M}^p_{(n,g,k+1)}$.

Moreover, we observe that different (up to isomorphism) admissible choices for $(C_1,{\bf g_1})$, $(C_2,{\bf g_2})$ correspond to the same component $\mathcal{M}^p_{(n,g,k+1)}$. In \verb|Table 1| we give one example for each component. It is remarkable that one can always construct this example choosing $C_2\simeq D_p$, ${\bf g_2}\simeq \delta_p$. In \verb|Table 1| one can find: the properties which characterize the pair $(C_1,{\bf g_1})$; the singularities of $(C_1\times D_p)/({\bf g_1}\times \delta_p)$; the value of $K_X^2$; the fixed locus $(n,g,k+1)$ of the automorphism induced on $S$;  and this identifies $\mathcal{M}^p_{(n,g,k+1)}$ whose dimension is $m$.\endproof

The quotients $(C_1\times D_3)/({\bf g_1}\times \delta_3)$ which admit a minimal model which is a K3 surface are classified in \cite[Remark 3.1]{GvG} and are listed in \verb|Table 1|.

\begin{eqnarray*}
\begin{array}{|c|c|c|c|c|c|c|c|}
\hline
p&g(C_1)&(\alpha_1,\ldots \alpha_{p-1})&(a_1,\ldots, a_{p-1})& Sing(C_1\times D_p/{\bf g_1}\times \delta_p)&K_X^2&(n,g,k+1)&m\\
\hline
3&4&(3,1)&(0,6)&\left(\frac{1}{3}\right)^{18}&-6&(6,0,3)&3\\
\hline
3&3&(2,1)&(1,4)&\left(\frac{1}{3}\right)^{12}, \left(\frac{2}{3}\right)^{3} &-4&(7,0,4)&2\\
\hline
3&2&(1,1)&(2,2)&\left(\frac{1}{3}\right)^{6}, \left(\frac{2}{3}\right)^{6} &-2&(8,0,5)&1\\
\hline
3&1&(0,1)&(3,0)&\left(\frac{2}{3}\right)^{9} &0&(9,0,6)&0\\
\hline 
5&6&(3,2,1,0)&(0,0,0,5)&\left(\frac{1}{5}\right)^{10}, \left(\frac{3}{5}\right)^5&-12&(7,0,1)&2\\
\hline
5&4&(2,1,1,0)&(0,1,0,3)&\left(\frac{1}{5}\right)^{6}, \left(\frac{2}{5}\right)^5, \left(\frac{4}{5}\right)&-8&(10,0,2)&1\\
\hline
5&2&(1,0,1,0)&(0,2,0,1)&\left(\frac{1}{5}\right)^{2}, \left(\frac{2}{5}\right)^5, \left(\frac{4}{5}\right)^2&-4&(13,0,3)&0\\
\hline
7&6&(2,2,1,1,0,0)&(0,0,0,0,1,3)&\left(\frac{1}{7}\right)^{6}, \left(\frac{3}{7}\right)^3, \left(\frac{4}{7}\right)^2, \left(\frac{5}{7}\right) &-14&(8,0,1)&1\\
\hline
7&3&(0, 1, 1, 0, 0, 1)&(0, 0, 0, 2, 1, 0 )&\left(\frac{1}{7}\right)^{2}, \left(\frac{2}{7}\right)^2, \left(\frac{3}{7}\right)^5  &-7&(13,0,2)&0\\
\hline
11&5&(1,1,0,1,1,&(0,0,1,0,0,&\left(\frac{1}{11}\right)^{2}, \left(\frac{2}{11}\right), \left(\frac{3}{11}\right)^2, &-13&(11,0,1)&0\\
&&0,0,1,0,0)&0,1,0,0,1)&
\left(\frac{4}{11}\right), \left(\frac{5}{11}\right), \left(\frac{7}{11}\right)^2
 &&&\\
\hline
13&6&(1,1,1,0,1,1,&(0,0,0,1,0,0,&\left(\frac{1}{13}\right)^{2}, \left(\frac{2}{13}\right), \left(\frac{3}{13}\right)^2, &-17&(9,0,1)&0\\
&&0,0,1,0,0,0)&0,0,0,1,0,1)&
\left(\frac{5}{13}\right), \left(\frac{6}{13}\right), \left(\frac{9}{13}\right)^2
 &&&\\
\hline
17&8&(1,1,1,1,1,1,0,1,&(0,0,0,0,0,0,0,0,&\left(\frac{1}{17}\right)^{2}, \left(\frac{4}{17}\right), \left(\frac{5}{17}\right), &-23&(7,-,0)&0\\
&&0,0,1,0,0,0,0,0,0)&0,0,0,0,1,0,1,1)&
\left(\frac{7}{17}\right)^2, \left(\frac{8}{17}\right), \left(\frac{9}{17}\right)^2
 &&&\\
 \hline
19&9&(1,1,1,1,1,1,&(0,0,0,0,0,0,&\left(\frac{1}{19}\right)^{2}, \left(\frac{3}{19}\right), \left(\frac{5}{19}\right)^2, &-25&(5,-,0)&0\\
&&0,1,1,0,0,1,)&0,0,0,0,0,0)&
\left(\frac{7}{19}\right), \left(\frac{9}{19}\right), \left(\frac{13}{19}\right)^2
 &&&\\
&&0,0,0,0,0,0)&0,0,1,1,0,1)&
&&&\\
\hline
\end{array}
\end{eqnarray*}
\begin{center}
\verb|Table 1|
\end{center}
\begin{rem} The example listed in \verb|Table 1| are all the examples obtained as described if $p\leq 7$. For $p\geq 11$ there are other admissible pairs $(C_1,{\bf g_1})$, $(C_2,{\bf g_2})$, such that $C_i\not\simeq D_p$, $i=1,2$, which correspond to the components $\mathcal{M}^{11}_{(11,0,1)}$, $\mathcal{M}^{13}_{(9,0,1)}$, $\mathcal{M}^{17}_{(7,-,0)}$, $\mathcal{M}^{19}_{(5,-,0)}$.\end{rem}

\begin{example}\label{subsec: example p=5}{\rm
As example (line 5 \verb|Table 1|) we consider the pairs $(C_1, {\bf g_1})$ and $(C_2, {\bf g_2})$ such that $|{\bf g_i}|=5$ and:\\
$\bullet$ $g(C_1)=6$, $C_1\ra C_1/ \langle {\bf g_1} \rangle$ is branched in 5 points and the local action of ${\bf g_1}$ near these points is $\zeta_5^{4}$. The dimension of the eigenspaces for the induced action in cohomology is $(\alpha_1,\ldots \alpha_4)=(3,2,1,0)$;\\
$\bullet$ $C_2\simeq D_5$, ${\bf g_2}=\delta_5$ (cf. Example \ref{ex: curve Dp}).\\
We will denote by $P_i$, $i=1,2,3,4,5$ the branch points of $C_1\ra \mathbb{P}^1$, by $Q_j$, $j=1,2,3$ the branch points of $C_2\ra \mathbb{P}^1$ and we assume the local action of ${\bf g_2}$ is the same near the points $Q_1$ and $Q_2$. \\
The singularities of the quotient $(C_1\times C_2)/({\bf g_1}\times {\bf g_2})$ are 10 singularities of type $\frac{1}{5}(1,1)$ (over $P_i\times Q_j$, $i=1\ldots,5$, $j=1,2$) and 5 of types $\frac{1}{5}(1,3)$ (over $P_i\times Q_3$, $i=1,\ldots,5$). The resolution of the 10 singularities of type $\frac{1}{5}(1,1)$ consists of the 10 $(-5)$-curves $\widetilde{B_{i,j}}$, $i=1,\ldots, 5$, $j=1,2$. The resolution of each singularity of type $\frac{1}{5}(1,3)$ consists of 2 curves, $\widetilde{B_{i,3}^h}$, $h=1,2$, meeting in one point and with self intersection $-2$ and $-3$, respectively. The resolution of the singularities is as in Figure 1.
\begin{center}
\includegraphics[totalheight=7 cm]{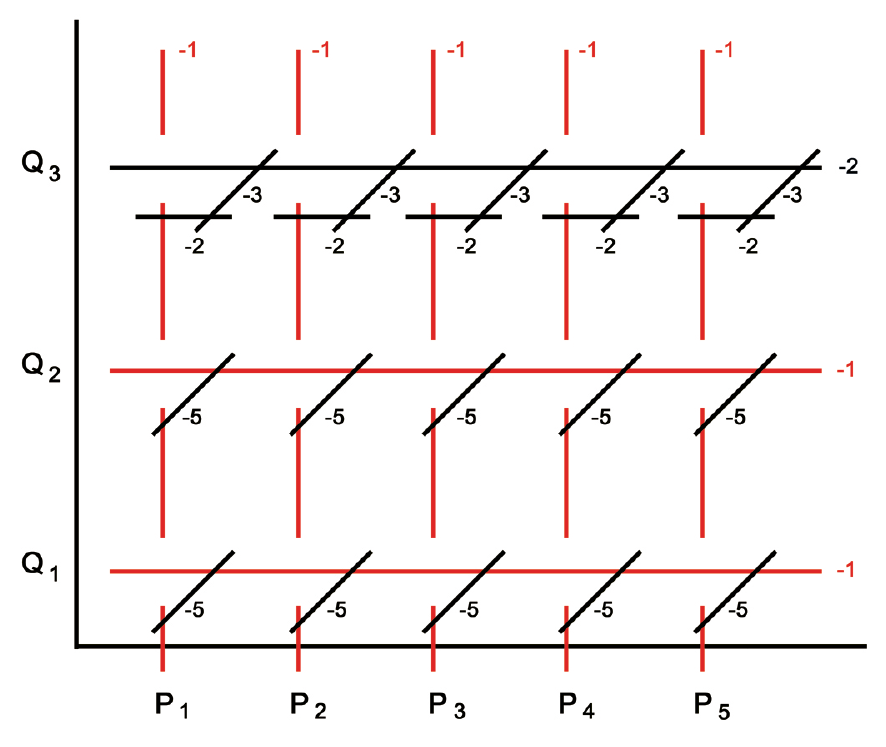} \\
Figure 1.
\end{center}
The central components $F_{P_i}$ of the 5 reducible fibers of the fibration $\pi_1\colon X\ra C_1/\langle {\bf g_1} \rangle$ are $(-1)$-curves. The central components $F_{Q_j}$ of the reducible fibers of $\pi_2\colon X\ra C_2/\langle {\bf g_2} \rangle$ over $Q_j$ with $j=1,2$ are $(-1)$-curves and the central component $F_{Q_3}$ of $\pi_2$ over $Q_3$ is a $(-2)$-curve (cf. Proposition \ref{prop: selfinter central}).\\

In order to construct the minimal model $S$ we first contract the 5 $(-1)$-curves $F_{P_i}$ (we call this contraction map $\sigma_1$). The images $B_{i,3}^1$ of the curves $\widetilde{B_{i,3}^1}$ are 5 $(-1)$-curves and thus we contract them (we call this contraction map $\sigma_2$). After the contraction of the 2 $(-1)$-curves $F_{Q_1}$ and $F_{Q_2}$ we are left with only $(-2)$-curves from the configuration we started with. We call this surface $S$ (see Figure 2).
\begin{center}
\includegraphics[totalheight=7 cm]{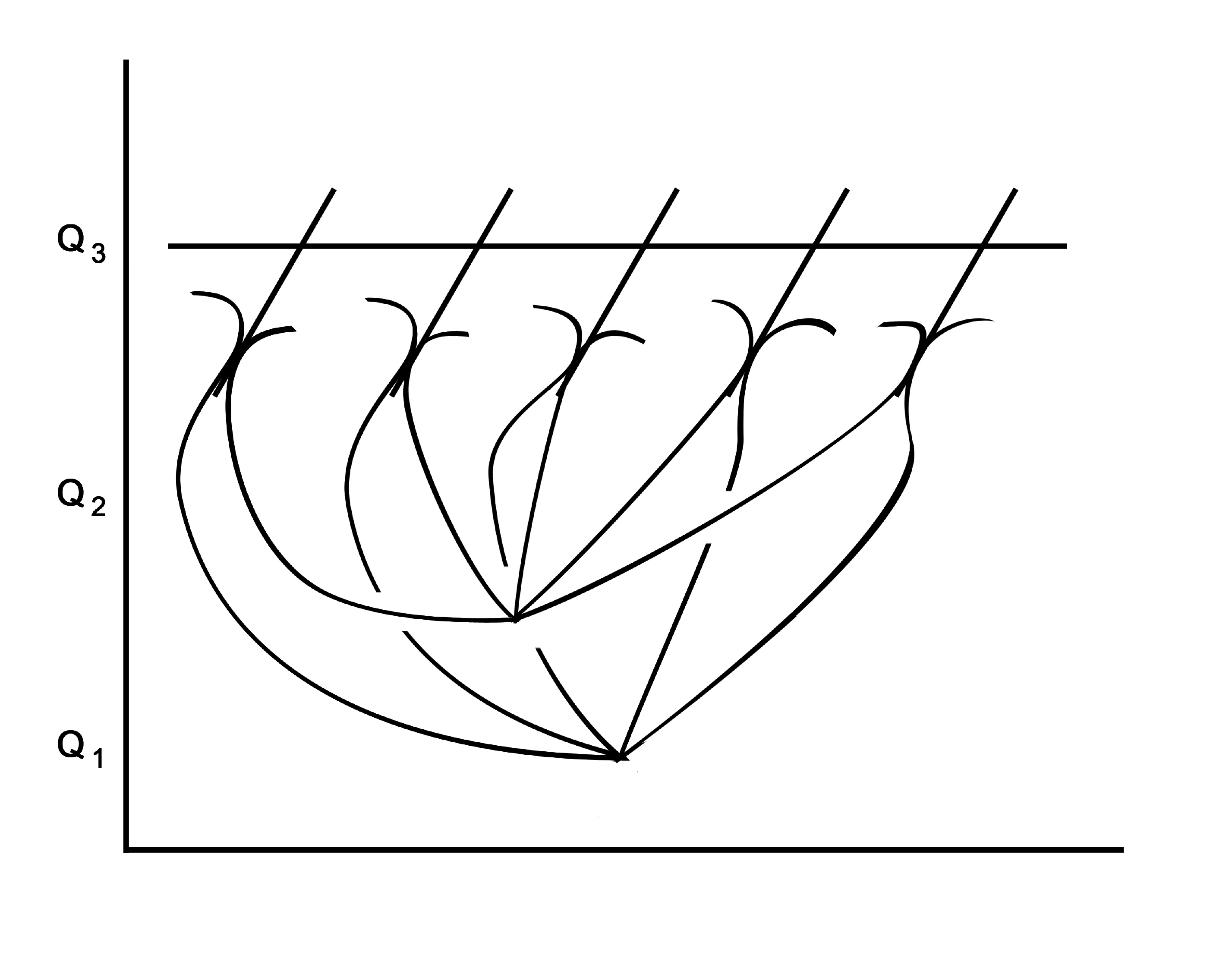} \\
Figure 2.
\end{center}
 We contracted 12 curves and since the canonical bundle of $X$ has self-intersection $K_X^2=-12$, we obtain $K_S^2=0$. Now we verify the surface $S$ satisfy the hypothesis of Lemma \ref{lem: trivial canonical bundle} and this proves $S$ is a K3 surface:\\
By adjunction, $2g(C_1)-2=(K_X-F_2)F_2=K_XF_2=10$ and $2g(C_2)-2=(K_X-F_1)F_1=K_XF_1=2$.
Since $\sigma_1^*(B_{i,3}^1)=\widetilde{B_{i,3}^1}+F_{P_i}$, the exceptional divisor $E$ of the blow up $X\ra S$ is $$E=\sum_{i=1}^5(2F_{P_i}+\widetilde{B_{i,3}^1})+F_{Q_1}+F_{Q_2}.$$ The curves $F_{P_i}$ are sections of the fibration $\pi_2$ and $F_{Q_j}$ are sections of the fibration $\pi_1$, hence $F_{P_i}F_2=F_{Q_j}F_1=1$. The curves $\widetilde{B_{i,j}^h}$ are components of the reducible fibers of both the fibration $\pi_1$ and $\pi_2$, hence $\widetilde{B_{i,j}^h}F_1=\widetilde{B_{i,j}^h}F_2=0$. So $(K_X-E)F_i=0$ for $i=1,2$.\\

By Remark \ref{rem: induced automorphism} the fixed locus of $g_S$ on $S$ consists of one rational curve (the image of the central components $F_{Q_3}$) and 7 points (see Figure 2).

Since $C_2$ is a rigid curve and $C_1$ varies in an irreducible 2-dimensional family, we in fact construct an irreducible 2-dimensional family $\mathcal{F}$ of K3 surfaces $S$ with a non-symplectic automorphism of order 5 and fixed locus $(n,g,k+1)=(7,0,1)$. Thus $\mathcal{F}\subset\mathcal{M}^5_{(7,0,1)}$. By \cite{AST}, $\mathcal{M}^5_{(7,0,1)}$ is an irreducible 2-dimensional family of K3 surfaces with the required automorphism and hence $\mathcal{F}$ coincides $\mathcal{M}^5_{(7,0,1)}$.
}\end{example}


\section{K3 surfaces which are minimal models of $(C_1\times C_2)/(\Z/2p\Z)$}\label{sec: PQ with Z2p}

In Section \ref{sec: PQ with Zp} we described K3 surfaces which are minimal models of a product-quotient surfaces with group $\Z/p\Z$. If $p\leq 11$ we never get the maximal irreducible components of the moduli space of K3 surfaces with a non-symplectic automorphism of order $p$.  In order to find at least one maximal irreducible component of such moduli space, we consider product-quotient with the $\Z/2p\Z$.
\begin{theo}\label{theo: PQ with Z2p}
Let $p=3$ (resp. $3<p\leq13$, $13< p\leq 19$). The K3 surface $S$ admitting a non-symplectic automorphism of order $p$ which fixes at least 2 (resp. 1, 0) curves are all minimal models of product-quotient surfaces with the group $\Z/2p\Z$. In particular we obtain the irreducible maximal component  $\mathcal{M}^3_{(n,g,2)}$ (resp. $\mathcal{M}^p_{(n,g,1)}$, $\mathcal{M}^p_{(n,-,0)}$).\\ Moreover, we can always choose $(C_2,{\bf g_2})\simeq (D_p,\tau_p)$ and hence the non-symplectic automorphism on $S$ is always induced by $\id\times\delta_p$.\\
In \verb|Table 2| we list an example for each family assuming $(C_2,{\bf g_2})\simeq (D_p,\tau_p)$.
Viceversa all the K3 surfaces which are minimal models of $(C_1\times C_2)/({\bf g_1}\times {\bf g_2})$, $|{\bf g_1}\times {\bf g_2}|=2p$ admit a non-symplectic automorphism of order $p$ whose fixed locus is one of those listed in \verb|Table 2|.
\end{theo}
\proof The proof is analogous to the one of Theorem \ref{theo: PQ with Zp}. We remark that in \verb|Table 2| one can find: the properties which characterize the pair $(C_1,{\bf g_1})$; the singularities of $(C_1\times D_p)/({\bf g_1}\times \tau_p)$; the value of $K_X^2$; the fixed locus $(n,g,k+1)$ of the automorphism induced on $S$;  and this identifies $\mathcal{M}^p_{(n,g,k+1)}$ whose dimension is $m$.\endproof

\begin{cor}\label{cor: autom 2p} All the K3 surfaces admitting a non-symplectic automorphism of order $p=3$ (resp. $3<p\leq 13$, $13< p\leq 19$) which fixes at least 2 (resp. 1,0) curves admit a non-symplectic automorphism of order $2p$ induced by $\id\times\tau_p$.\end{cor}
We observe that the results in Corollary \ref{cor: autom 2p} were already proved in \cite{Dillies} for $p=3$ and \cite{GS} for $p>3$.

{\scriptsize
\begin{eqnarray*}
\begin{array}{|c|c|c|c|c|c|c|c|}
\hline
p&g(C_1)&(\alpha_1,\ldots \alpha_{2p-1})&(a_1,\ldots, a_{2p-1})& Sing(C_1\times D_p/{\bf g_1}\times \tau_p)&K_X^2&(n,g,k+1)&m\\
\hline
3&25&(9,7,5,3,1)&(0,12,0,0,0)&\left(\frac{1}{6}\right)^{12}, \left(\frac{1}{3}\right)^{12},
\left(\frac{1}{2}\right)^{12} 
 &-36&(0,5,2)&9\\
\hline
3&22&(8,6,4,3,1)&(0,10,2,0,0)&\left(\frac{1}{6}\right)^{10},
\left(\frac{1}{3}\right)^{13},
\left(\frac{1}{2}\right)^{10} 
 &-31&(1,4,2)&8\\
\hline 3&19&(7,5,3,3,1)&(0,8,4,0,0)&\left(\frac{1}{6}\right)^{8},
\left(\frac{1}{3}\right)^{14},
\left(\frac{1}{2}\right)^{8} 
 &-26&(2,3,2)&7\\
\hline 3&16&(6,4,2,3,1)&(0,6,6,0,0)&\left(\frac{1}{6}\right)^{6},
\left(\frac{1}{3}\right)^{15},
\left(\frac{1}{2}\right)^{6} 
 &-21&(3,2,2)&6\\
\hline 3&17&(6,5,3,2,1)&(0,8,0,2,0)&\left(\frac{1}{6}\right)^{8},
\left(\frac{1}{3}\right)^{8},\left(\frac{2}{3}\right)^{3},
\left(\frac{1}{2}\right)^{8} 
 &-24&(3,3,3)&6\\
\hline 3&13&(5,3,1,3,1)&(0,4,8,0,0)&\left(\frac{1}{6}\right)^{4},
\left(\frac{1}{3}\right)^{16},
\left(\frac{1}{2}\right)^{4} 
 &-16&(4,1,2)&5\\
\hline 3&14&(5,4,2,2,1)&(0,6,2,2,0)&\left(\frac{1}{6}\right)^{6},
\left(\frac{1}{3}\right)^{9},\left(\frac{2}{3}\right)^3,
\left(\frac{1}{2}\right)^{6} 
 &-19&(4,2,3)&5\\
\hline 3&15&(5,4,3,2,1)&(1,7,0,0,0)&\left(\frac{1}{6}\right)^{7},\left(\frac{5}{6}\right),
\left(\frac{1}{3}\right)^{7},\left(\frac{2}{3}\right),
\left(\frac{1}{2}\right)^{8} 
 &-21&(4,3,4)&5\\
\hline 3&10&(4,2,0,3,1)&(0,2,10,0,0)&\left(\frac{1}{6}\right)^{2},
\left(\frac{1}{3}\right)^{17},
\left(\frac{1}{2}\right)^{2} 
 &-11&(5,0,2)&4\\
\hline 3&11&(4,3,1,2,1)&(0,4,4,2,0)&\left(\frac{1}{6}\right)^{4},
\left(\frac{1}{3}\right)^{10},\left(\frac{2}{3}\right)^3,
\left(\frac{1}{2}\right)^{4} 
 &-14&(5,1,3)&4\\
\hline 3&12&(4,3,2,2,1)&(1,5,2,0,0)&\left(\frac{1}{6}\right)^{5},\left(\frac{5}{6}\right),
\left(\frac{1}{3}\right)^{8},\left(\frac{2}{3}\right),
\left(\frac{1}{2}\right)^{6} 
 &-16&(5,2,4)&4\\
 \hline 3&7&(3,1,0,2,1)&(0,1,8,0,3)&\left(\frac{1}{6}\right),
\left(\frac{1}{3}\right)^{13},
\left(\frac{1}{2}\right)^{5} 
 &-7&(6,0,3)&3\\
 \hline 3&9&(3,2,1,2,1)&(1,3,4,0,0)&\left(\frac{1}{6}\right)^3,\left(\frac{5}{6}\right),
\left(\frac{1}{3}\right)^{9},\left(\frac{2}{3}\right),
\left(\frac{1}{2}\right)^{4} 
 &-11&(6,1,4)&3\\
 \hline 3&4&(2,0,0,1,1)&(0,0,6,0,6)&
\left(\frac{1}{3}\right)^{9},
\left(\frac{1}{2}\right)^{8} 
 &-3&(7,0,4)&2\\
 \hline 3&7&(2,2,1,1,1)&(1,3,0,2,0)&
\left(\frac{1}{6}\right)^3,\left(\frac{5}{6}\right),
\left(\frac{1}{3}\right)^3,\left(\frac{2}{3}\right)^4,
\left(\frac{1}{2}\right)^{4}  
 &-9&(7,1,5)&2\\
 \hline 3&3&(1,0,0,1,1)&(1,0,4,0,3)&\left(\frac{5}{6}\right),
\left(\frac{1}{3}\right)^{6},\left(\frac{2}{3}\right),
\left(\frac{1}{2}\right)^{5} 
 &-2&(8,0,5)&1\\
 \hline 3&5&(1,1,1,1,1)&(2,2,0,0,0)&\left(\frac{1}{6}\right)^2,\left(\frac{5}{6}\right)^2,
\left(\frac{1}{3}\right)^2,\left(\frac{2}{3}\right)^2,
\left(\frac{1}{2}\right)^{4} 
 &-6&(8,1,6)&1\\
 \hline 3&1&(0,0,0,0,1)&(0,1,2,0,3)&\left(\frac{5}{6}\right),
\left(\frac{2}{3}\right)^{4},
\left(\frac{1}{2}\right)^{5} 
 &0&(9,0,6)&0\\
\hline
5&22&(4,1,0,2,&(0,0,6,0,&\left(\frac{1}{10}\right)^6, \left(\frac{1}{5}\right)^2,
 &-28&(1,2,1)&4\\
&&2,4,5,3,1)&0,0,2,0,0)&
 \left(\frac{2}{5}\right)^6, \left(\frac{3}{5}\right),\left(\frac{1}{2}\right)^6
 &&&\\
 \hline
 5&17&(3,2,0,2,&(0,0,4,0,&\left(\frac{1}{10}\right)^4, \left(\frac{1}{5}\right)^3, \left(\frac{2}{5}\right)
 &-21&(4,1,1)&3\\
&&1,1,4,3,1)&2,0,2,0,0)&
\left(\frac{3}{5}\right)^6, \left(\frac{1}{2}\right)^4
 &&&\\
    \hline
 5&12&(3,1,1,2,&(0,0,1,1,&\left(\frac{1}{10}\right), \left(\frac{7}{10}\right),\left(\frac{1}{5}\right)^7
 &-13&(7,1,2)&2\\
&&0,2,2,1,0)&0,0,6,0,0)&
\left(\frac{2}{5}\right)^3, \left(\frac{3}{5}\right),
\left(\frac{4}{5}\right), \left(\frac{1}{2}\right)^2 &&&\\
   \hline
 5&8&(2,1,1,1,&(0,0,1,1,&\left(\frac{1}{10}\right), \left(\frac{7}{10}\right),\left(\frac{1}{5}\right)^3
 &-9&(10,0,2)&1\\
&&0,1,1,1,0)&0,0,2,2,0)&
\left(\frac{2}{5}\right)^3, \left(\frac{3}{5}\right),
\left(\frac{4}{5}\right),
 &&&\\
 &&&& \left(\frac{1}{2}\right)^2 &&&\\
  \hline
 5&4&(1,0,1,1,&(0,1,0,1,&\left(\frac{3}{10}\right), \left(\frac{9}{10}\right),
 &-5&(13,0,3)&0\\
&&0,1,0,0,0)&0,0,2,0,0)&
\left(\frac{1}{5}\right)^3, \left(\frac{3}{5}\right)^2,
\left(\frac{1}{2}\right)^2
 &&&\\
\hline
7&19&(2,2,3,1,0,2,&(0,0,0,4,0,0,&\left(\frac{1}{14}\right)^4,
\left(\frac{4}{7}\right)^2,
 &-25&(3,1,1)&2\\
&&1,0,3,1,0,3,1)&0,0,2,0,0,0,0)&
\left(\frac{5}{7}\right)^5, 
\left(\frac{1}{2}\right)^4
 &&&\\
  \hline
7&12&(1,1,2,1,0,1,&(0,0,0,1,1,0,&\left(\frac{1}{14}\right),
\left(\frac{5}{14}\right),\left(\frac{1}{7}\right)^3,
 &-14&(8,1,2)&1\\
&&0,1,2,1,0,1,1)&2,0,0,0,2,0,0)&
\left(\frac{3}{7}\right),
\left(\frac{4}{7}\right),\left(\frac{5}{7}\right),^3
 &&&\\
 &&&&\left(\frac{1}{2}\right)^2 &&&\\
 \hline
7&6&(0,1,1,0,0,0,&(1,0,0,1,0,0,&\left(\frac{1}{14}\right),
\left(\frac{9}{14}\right),\left(\frac{1}{7}\right),
 &-10&(13,0,2)&0\\
&&0,0,1,1,0,1,1)&2,0,0,0,0,0,0)&
\left(\frac{3}{7}\right)^3, \left(\frac{6}{7}\right),
\left(\frac{1}{2}\right)^2
 &&&\\
 \hline
 11&21&(1,2,2,0,2,1,0,1,0,2&(0,0,0,0,0,3,0,0,1&\left(\frac{1}{22}\right)^3, \left(\frac{19}{22}\right),
 &-31&(2,1,1)&1\\
&&1,0,2,0,2,1,0,1,0,&0,0,0,0,0,0,0,0,0,0&
\left(\frac{6}{11}\right),
\left(\frac{9}{11}\right)^3, \left(\frac{1}{2}\right)^4
 &&&\\
&& 2,1) &0,0) &&&&\\
 \hline
 11&10&(1,0,1,1,0,1,0,0,0,1&(0,0,0,0,0,0,1,0,1,&\left(\frac{3}{22}\right),
 \left(\frac{7}{22}\right),\left(\frac{1}{11}\right)^2,
 &-13&(11,0,1)&0\\
&&0,1,1,0,1,1,1,0,0,&0,0,0,0,0,0,0,0,0,2,&
\left(\frac{2}{11}\right), \left(\frac{4}{11}\right),
\left(\frac{9}{11}\right),
 &&&\\
&& 0,0) &0,0) &\left(\frac{1}{2}\right)^2&&&\\
 \hline
 13&12&(0,1,1,0,1,0,0,1,0,1&(0,0,0,0,0,0,1,0,1&\left(\frac{1}{26}\right), \left(\frac{5}{26}\right),\left(\frac{3}{13}\right)^2,
 &-18&(9,0,1)&0\\
&&0,0,0,0,1,1,1,0,1,&0,0,0,0,0,0,0,0,0,2&
\left(\frac{4}{13}\right), \left(\frac{6}{13}\right),
\left(\frac{7}{13}\right),
 &&&\\
&& 0,0,1,0,1,1,1) &0,0,0,0,0,0) &\left(\frac{1}{2}\right)^2&&&\\
 \hline
 17&16&(1,1,1,0,1,0,1,0,0,1&(0,0,0,0,0,0,0,0,1&\left(\frac{1}{34}\right), \left(\frac{23}{34}\right),
 &-22&(7,-,0)&0\\
&&0,1,0,0,1,1,0,0,0,&1,0,0,0,0,0,0,0,0,0&
\left(\frac{5}{17}\right)^2, \left(\frac{7}{17}\right)^2,
\left(\frac{15}{17}\right),
 &&&\\
&& 0,1,1,1,0,1,0,0) &0,0,0,2,0,0,0,0) &\left(\frac{1}{2}\right)^2&&&\\
&& 1,0,1,0,1,0) &0,0,0,0,0,0) &&&&\\
 \hline
 19&18&(1,0,1,1,1,1,0,0,1,1&(0,0,0,0,0,0,0,0,0&\left(\frac{3}{38}\right), \left(\frac{11}{38}\right),\left(\frac{1}{19}\right)^2,
 &-22&(5,-,0)&0\\
&&0,1,0,0,0,0,0,1,0,&1,0,0,0,1,0,0,0,0,0&
 \left(\frac{4}{19}\right),
\left(\frac{8}{19}\right), \left(\frac{17}{19}\right)
 &&&\\
&& 1,1,0,1,1,1,1,1,0) &0,0,0,0,0,0,0,0,0) &\left(\frac{1}{2}\right)^2&&&\\
&& 0,0,1,1,0,0,0,0,0) &0,0,0,0,2,0,0) &&&&\\
\hline
\end{array}
\end{eqnarray*}
\begin{center}
\verb|Table 2|
\end{center}}

\begin{rem}{\rm The K3 surfaces constructed in the \verb|Table 2| admit an automorphism of order $2p$, induced by $\id\times \tau_p$ (Corollary \ref{cor: autom 2p}) and one of order 2, induced by $\id\times\iota_p$. The fixed locus of these automorphisms can be computed case by case. In Section \ref{sec: equation} we compute it in certain cases, by using a projective model of the surfaces.}\end{rem}

\subsection{Intermediate quotients}
The $2p:1$ map $C_1\times D_p\ra (C_1\times D_p)/{(\bf g_1\times \tau_p)}$ clearly factorizes through $$C_1\times D_p\stackrel{p:1}{\longrightarrow} (C_1\times D_p)/({\bf g_1}\times \tau_p)^2\stackrel {2:1}{\longrightarrow}(C_1\times D_p)/({\bf g_1}\times \tau_p).$$ This induces a 2:1 rational map between the minimal model, $Q$, of $(C_1\times D_p)/({\bf g_1}\times \tau_p)^2$ and the K3 surface $S$. In particular $Q$ is a 2-cover of a K3 surface. We observe that $p_g(Q)\geq p_g(S)$. This immediately implies that the Kodaira dimension of $Q$, $k(Q)$, is non negative. The following examples show that all the other three possibilities, $k(Q)=0,1,2$, appear in our classification.
First we notice that the genus of the quotient $C_1/{\bf g_1}^2$ is $\alpha_p$ and so $q(Q)=\alpha_p$.

Let us consider the line 18 of \verb|Table 2|, it corresponds to the quotient $(D_3\times D_3)/(\tau_3^5\times \tau_3)$. The quotient $(D_3\times D_3)/(\tau_3^5\times \tau_3)^2$ is isomorphic to $(D_3\times D_3)/(\delta_3\times \delta_3^2)$. The minimal model of such a surface is described in \cite{SI} (see Example \ref{example: quotient hodge numbers 0,i} the fourth line of \verb|Table 1|) and is a K3 surface. In particular in this case $k(Q)=0$.

Let us consider the line 9 of \verb|Table 2|. The map $C_1\ra C_1/{\bf g_1}^2\simeq \mathbb{P}^1$ is branched in 12 points and an equation of $C_1$ is $w^3=p_{12}(t)$ where $p_{12}(t)$ is a polynomial with 12 simple roots. With the same method we will apply in Section \ref{sec: equation}, case $p=3$, we obtain the equation $y^2=x^3+p_{12}^2(t)$ of $(C_1\times D_p)/({\bf g_1}\times \tau_p)^2$.  So the surface $Q$ admits an elliptic fibration, its birational invariant are $q(Q)=0$, $p_g(Q)=3=\alpha_2+\alpha_5$, and we obtain $k(Q)=1$.

Let us consider the line 19 of \verb|Table 2|. It corresponds to the quotient $(C_1\times D_5)/({\bf g_1}\times \tau_5)$ where $g(C_1)=22$. Let $Y$ be the minimal resolution of $(C_1\times D_5)/({\bf g_1}\times \tau_5)^2$. As in Example \ref{subsec: example p=5} , one proves that the singularities of $(C_1\times D_5)/({\bf g_1}\times \tau_5)^2$ are 10 singularities of type $\frac{1}{5}(1,1)$, 12 singularities of type $\frac{1}{5}(1,2)$ and 2 singularities of type $\frac{1}{5}(1,3)$. The computation of $K_Y^2$ can be done as explained in Remark \ref{rem: PQ properties} (4) and it gives $K_Y^2=10>0$. Since $Q$ is the minimal model of $Y$, $K_Q^2\geq K_Y^2$ and we conclude that $K_Q^2>0$, $Q$ is a surface of general type, and so $k(Q)=2$.


\section{Equations}\label{sec: equation}


\subsection{Automorphisms of order $p=3$} In \cite[Proposition
4.2]{AS08} it is proved that every K3 surfaces admitting a
non-symplectic automorphism of order 3 fixing at least two curves
is in fact an isotrivial elliptic fibration with generic fiber
isomorphic to the elliptic curve $E_{\zeta_3}$ with complex
multiplication of order 3. Indeed, every such a K3 surface is
described as en elliptic K3 surface with an equation of type
$y^2=x^3+f_{12}(t)$. In view of our construction this is very
natural: we proved that every such a K3 surface is the minimal
model of the quotient $(C_1\times D_3)/({\bf g_1}\times \tau_3)$
where $D_3\simeq E_{\zeta_3}$ and $\tau_3$ are described in Example \ref{ex: curve
Dp} and $(C_1,{\bf g_1})$ varies. The maximal component is
obtained by $(C_1,{\bf g_1})$ as in the first line of the
\verb|Table 2|. In this case $C_1$ is a $6:1$ cover of
$\mathbb{P}^1$ whose ramification consists of 12 points of order
6. An equation of $C_1$ is $w^6=f_{12}(t)$, where $\deg(f_{12}(t))=12$ and $f_{12}(t)$ does not have multiple roots.
The local action near the fixed points is $-\zeta_3^2$ (see
\verb|Table 2|) and thus we can assume that the automorphism ${\bf
g_1}$ is ${\bf g_1}\colon (w,t)\mapsto (-\zeta_3^2w,t)$. The new
functions $x:=uw^2$, $y:=vw^3$ and $t$ are invariant for ${\bf
g_1}\times \tau_3$ and satisfy the equation
$$y^2=x^3+f_{12}(t).$$ Moreover, if $W$ is
the surface defined by this equation, then the generic fiber of
the map $C_1\times D_3\ra W$ consist of 6 points, thus we have the
following commutative diagram:
$$\begin{array}{rcl}
&C_1\times D_3&\\
6:1\swarrow&&\searrow 6:1\\
W\ \ \ &\dashrightarrow&(C_1\times D_3)/({\bf g_1}\times \tau_3)
\end{array}$$
which shows that $W$ and $S$ are birational and so $W$ is a
singular model of the K3 surface $S$. This construction was also considered in \cite[Example 3.11]{vG1}.

More in general the curve $C_1$ has an equation of type $w^6=f_{12}(t)$ where $f_{12}$ does
not admit roots with multiplicity greater then 5 and there exists
no a polynomial $h(t)$ such that either $f_{12}(t)=h^2(t)$ or
$f_{12}(t)=h^3(t)$.  If some of the roots
of $f_{12}(t)$ have multiplicity higher then 1, then the fixed
locus of $\id\times \delta_3$ changes and we obtain a member of a
more special family (cf. lines from 2 to 18 \verb|Table 2|).

We saw in Section \ref{sec: PQ with Zp} that certain K3 surfaces
admitting a non-symplectic automorphism of order 3, can be
obtained as quotient of $(C_1\times D_3)$ by an automorphism of
order 3, ${\bf g_1}\times\delta_3$. So we obtain a different
equation for these K3 surfaces. Indeed, in this case one can assume
$C_1$ to have the following equation $w^3=f_6(t)$, such that
$f_6(t)$ does not admit roots with multiplicity greater then 2. In
a very similar way as before this gives the following equation for
the quotient surface:
$$y^2=x^3+f_6(t)^2.$$ These equations were already considered in
\cite{GvG} and \cite{GS}.

\subsection{Automorphisms of order $p=5$}

In Theorem \ref{theo: PQ with Z2p} and \verb|Table 2| we proved
that the K3 surfaces admitting a non-symplectic
automorphism of order $5$ with at least one curve in the fixed locus are
the minimal models of quotients $(C_1\times D_5)/({\bf g_1}\times \tau_5)$
for a certain choices of the pair $(C_1,{\bf g_1})$. In particular the maximal component
(with fixed locus $(n,g(C),k+1)=(1,2,1)$) is obtained choosing $C_1$ to be a $10:1$
cover of $\mathbb{P}^1$ branched along 6 points of order $10$ and $1$ point of order $5$.
An equation of $C_1$ is $w^{10}=f_6(t)$ where $\deg(f_6(t))=6$ and $f_6(t)$ does not have multiple roots (we are assuming the branch point of order 5 is at infinity). The local action near the fixed points is $-\zeta_5^3$ (see \verb|Table 2|)
and so we can assume that the automorphism ${\bf g_1}$ is ${\bf g_1}:(w,t)\mapsto (-\zeta_5^2w,t)$.
The functions $x:=uw^2$, $y:=vw^5$ and $t$ are invariant under ${\bf g_1}\times \tau_5$ and satisfy the equation
$$y^2=x^5+f_6(t).$$
As in the case $p=3$, one shows that this gives in fact a (singular) model of the K3 surface $S$.
The equation exhibit $S$ as double cover of $\mathbb{P}^2_{[x,t]}$ branched along the (non homogenous) sextic $x^5+f_6(t)=0$.
A similar model for this K3 surface is described in \cite[Example 5.1]{AST}, where the relation with the curves $C_1$ and $D_5$ was not observed.\\
More in general, we observe that every curve $C_1$ in the
\verb|Table 2| admits an equation of the type $w^{10}=f_6(t)$ with
$f_6(t)$ which is not a square, such that ${\bf g_1}:(w,t)\ra
(-\zeta_5^2w,t)$. If $f_6(t)$ is generic we find the previous equation and
so the maximal component of the moduli space of K3 surfaces
admitting a non-symplectic automorphism of order 5 fixing at least
one curve. Specializations of  the polynomial $f_6(t)$ induce
specializations of the K3 surface $S$. For example the line 20 of
the \verb|Table 2| corresponds to the curve $C_1$ given by
$w^{10}=t^2g_4(t)$, $\deg(g_4(t))=4$ and $g_4$ does not have multiple roots. The corresponding K3 surface is
to the double cover of $\mathbb{P}^2_{[x,t]}$ branched along the
sextic
$x^5+t^2g_4(t)=0$, which has a singular point of type $A_4$ in the point $(x,t)=(0,0)$.\\
The automorphism $\tau_5$ (resp. $\delta_5$) on $D_5$ induces the non-symplectic
automorphism $\id\times \tau_5$ (resp.  $\id\times \delta_5$) of order 10 (resp. 5) on
the K3 surface $S$ which acts on the coordinates $(x,y,t)$ as $(x,y,t)\ra (\zeta_5 x,-y, t)$
(resp. $(x,y,t)\ra (\zeta_5 x,y, t)$). The fixed locus of $\id\times \delta_5$ consists of
one curve of genus 2 if $f_6(t)$ is generic and specializes to different fixed locus when $f_6(t)$ specializes
(see also \cite[Example 5.1]{AST}). \\
We observe that the non-symplectic automorphism $(id\times\tau_5)^5$ of order 2 is exactly the cover involution of the double cover of $\mathbb{P}^2$ and this allows one to compute easily its fixed locus.

We saw in Section \ref{sec: PQ with Zp} that certain K3 surfaces
admitting a non-symplectic automorphism of order $5$, can be
obtained from the quotient
$(C_1\times D_5)/ ({\bf g_1}\times\delta_5)$.
So we obtain a different equation for these K3 surfaces. The surfaces obtained in this way are listed in \verb|Table 1|.
In the case of the 1-dimensional  and 0-dimensional families
the equation of the curve $C_1$ is $w^5=f_3(t)$, $\deg(f_3(t))=3$ and $f_3$ is not a cube, and the automorphism
is ${\bf g_1}\colon(w,t)\ra (\zeta_5^2w,t)$. The functions $x:=uw^2$, $y:=vw^5$ and $t$ are invariant and give
a (singular) model of the K3 surface $S$, with equation $y^2=x^5s+f_3^2(t)$.\\

Every K3 surface that is the double cover of $\mathbb{P}^2$ branched along a sextic can be viewed as a hypersurface in the weighted projective space $W\mathbb{P}(3,1,1,1)$. In particular the homogeneous equation of $S$ can be written as $y^2=x^5s+f_6(t:s)$ where $(y:x:s:t)$ are the homogeneous coordinates of $W\mathbb{P}(3,1,1,1)$ ($y$ is the coordinate of weight 3).This remark will be useful in view of the equations we found in cases $p=7,11$.\\

\subsection{Automorphisms of order $p=7$}

In Theorem \ref{theo: PQ with Z2p} and \verb|Table 2| we proved that the K3 surfaces
admitting a non-symplectic automorphism of order 7 with at least
one curve in the fixed locus are the minimal model of the
quotient $(C_1\times D_7)/({\bf g_1}\times \tau_7)$ for a certain
choice of the pair $(C_1,{\bf g_1})$. In particular the maximal
component (with fixed locus $(n,g(C),k+1)=(3,1,1)$) is obtained
choosing $C_1$ to be a $14:1$ cover of $\mathbb{P}^1$ branched
along 4 points of order 14 and 1 point of order 7. An equation of
$C_1$ is $w^{14}=t(t-1)(t-\lambda_1)(t-\lambda_2)$. The local
action near the fixed points is $-\zeta_7^4$ (see \verb|Table 2|) and
thus the automorphism is ${\bf
g_1}:(w,t)\mapsto (-\zeta_7^3w,t)$. The functions $x:=uw^2$,
$y:=vw^7$, $t$ are invariant for ${\bf g_1}\times \tau_7$ and
satisfy the equation $y^2=x^7+t(t-1)(t-\lambda_1)(t-\lambda_2).$
As in case $p=3$, one shows that this gives in fact a (singular)
model of the K3 surface $S$: The equation can be homogeneized to
\begin{equation}\label{equation:
p=7}y^2=x^7s+t(t-s^2)(t-\lambda_1s^2)(t-\lambda_2s^2)\subset
W\mathbb{P}(4,2,1,1)_{(y:t:x:s)}.\end{equation}
In order to show that the equation \eqref{equation: p=7} corresponds in fact to a
(singular model of a) K3 surface we observe that the surface defined by \eqref{equation: p=7}
is well formed (cf. \cite[Definition 6.9]{Fl00}) and quasismooth (cf. \cite[Definition 6.3]{Fl00}).
If a hypersurface $Z$ of degree $d$ in a weighted projective space $W\mathbb{P}(a_0,a_1,a_2,a_3)$ is well formed and quasismooth,
then the adjunction formula generalizes and the canonical sheaf is $\omega_Z\simeq \mathcal{O}_Z(d-\sum_{i=0}^3 a_i)$ (cf. \cite[Paragraph 6.14]{Fl00}).
In particular if $d=\sum_{i=0}^3$, then $Z$ is a K3 surface and so the surface defined by \eqref{equation: p=7} is a singular model of a K3 surface.\\
We recall that the generic hypersurface of degree $8$ in $W\mathbb{P}(4,2,1,1)$ is a singular model of a K3 surface (\cite[Section 4.5]{R})
 with two singularities of type $\frac{1}{2}(1,1)$ at the points $(1:1:0:0)$, $(-1:1:0:0)$. The surface defined by \eqref{equation: p=7} (which is not general)
 has no other singular points.\\
The automorphism induced on $S$ by $\id\times \tau_7$ acts on the
coordinates of $W\mathbb{P}(4,2,1,1)$ in the following way:
$(y:t:x:s)\mapsto (-y:t:\zeta^2_7x:s)$. It has order 14 and its
fixed locus consists of 5 points: $(0:0:1:0)$, $(0:0:0:1)$,
$(0:1:0:1)$, $(0:\lambda_1:0:1)$, $(0:\lambda_2:0:1)$. The
singular points of $W\mathbb{P}(4,2,1,1)$ are switched by the
automorphism.
The fixed locus of the non-symplectic automorphism of order 7 induced by $\id\times \delta_7$ consists of
the point $(0:0:1:0)$ and of the curve $y^2=t(t-s^2)(t-\lambda_1s^2)(t-\lambda_2s^2)\subset W\mathbb{P}(4,2,1)$. The well formed expression (cf. \cite[Definition 5.11]{Fl00}) of this curve is $y^2=t(t-s)(t-\lambda_1 s)(t-\lambda_2s)\subset W\mathbb{P}(2,1,1)$ (cf. \cite[Lemma 5.7]{Fl00}) and the genus of such a curve is $1$ (e.g., \cite[Corollary 3.4.4]{D82}). \\
We observe that the singular points of $W\mathbb{P}(4,2,1,1)$ are contained in the fixed locus of the automorphism. By equation \eqref{equation: p=7} one sees that $S$ is a $2:1$ ramified cover of $W\PP(2,1,1)$, with branch locus given by $x^7s+t(t-s^2)(t-\lambda_1s^2)(t-\lambda_2s^2)=0$. The weighted projective plane $W\PP(2,1,1)$ has a natural embedding in $\PP^3$ with coordinates $(x_0:x_1:x_2:x_3)=(t:x^2:s^2:xs)$, whose image is a cone $Q$ of equation $x^2_3=x_1x_2$. The branch locus of the covering is now given by the intersection of $Q$ and the curve $x^3_1x_3+x_0(x_0-x_2)(x_0-\lambda_1x_2)(x_0-\lambda_2x_2)=0$, which does not pass through the vertex of $Q$. The automorphism descends to $\PP^3$ with the action $(x_0:x_1:x_2:x_3) \mapsto (x_0:\zeta^2_7 x_1: x_2:\zeta_7x_3)$. The fixed locus is the isolated point $(0:1:0:0)$ and the curve $x_1=x_3=0$, which passes through the vertex of $Q$.
We blow up the vertex of $Q$ introducing a copy of $\mathbb{P}^1$. The induced automorphism leaves invariant the exceptional divisor $E$ and fixes the strict transform of the fixed curve $B$. Since it restricts to an automorphism of $E$, it fixes two points on it, one of them is $E\cap B$. Above the other fixed point on $E$ we find two fixed point on $S$.

The fixed locus of the non-symplectic involution $\id\times \iota$ is the curve
$x^7s+t(t-s^2)(t-\lambda_1s^2)(t-\lambda_2s^2)\subset W\mathbb{P}(2,1,1)$. This is a curve of genus 9 in $W\mathbb{P}(2,1,1)$ by \cite[Corollary 3.4.4]{D82}.\\

In \verb|Table 1| we showed that certain K3 surfaces admitting a non-symplectic automorphism of order $7$, can be obtained from the quotient $(C_1\times D_7)/ ({\bf g_1}\times\delta_7)$.
In the case of the 0-dimensional family the equation of the curve $C_1$ is $w^7=t(t-1)$ and the automorphism is ${\bf g_1}\colon(w,t)\ra (\zeta_7^3w,t)$.  The functions $x:=uw^2$, $y:=vw^7$ and $t$ are invariant and gives a (singular) model of the K3 surface $S$, with equation $y^2=x^7s+t^2(t-s^2)^2\subset W\mathbb{P}(4,2,1,1)_{(y:t:x:s)}$.

\subsection{Automorphisms of order $p=11$}
If $p=11$, one can obtain an equation for a (singular model) of
$S$, minimal model of $(C_1\times D_{11})/({\bf g_1}\times
\tau_{11})$, as in cases $p=3,5,7$: an equation for the curve
$C_1$ is $w^{11}=t(t-1)(t-\lambda)$ (where if $\lambda\neq 0,1$
the curve $C_1$ is the one described in line 27 of  \verb|Table 2|, if
either $\lambda=1$ or $\lambda=0$, the curve $C_1$ is the one
described in line 28 of  \verb| Table 2|) and the automorphism ${\bf
g_1}\colon(w,t)\ra (-\zeta_{11}^5,t)$. An equation of $S$ is
$y^2=x^{11}s-t(t-s^4)(t-\lambda s^4)$ where $y:=vw^{11}$, $t$,
$x:=uw^2$ and $s$ are coordinates of  the weighted projective
space $W\mathbb{P}(6,4,1,1)$. As in case $p=7$, one shows that
this equation define in fact a singular model $W$ of a K3 surface.
The surface $W$ is singular in the point $(1:1:0:0)$.

The automorphism $\id\times\tau_{11}$ induces the non-symplectic automorphism $(y:t:x:s)\mapsto (-y:t:\zeta_{11}x:s)$ on the surface $W$ whose fixed locus consists of the points $(0:0:0:1)$, $(0:1:0:1)$, $(0:\lambda:0:1)$, $(0:0:1:0)$ (which are all distinct if $\lambda\neq 0$ and $\lambda\neq 1$).The point $(0:0:1:0)$ is a singular point of type of the surface.

The automorphism $\id\times\delta_{11}$ induces the non-symplectic automorphism $(y:t:x:s)\mapsto (y:t:\zeta_{11}x:s)$ on the surface $W$ whose fixed locus consists of the point $(0:0:1:0)$ and of the curve $y^2=t(t-s^4)(t-\lambda s^4)\subset W\mathbb{P}(6,4,1)$. The well formed expression of this curve is $y^2=t(t-s^2)(t-\lambda s^2)\subset W\mathbb{P}(3,2,1)$ which is quasismooth if $\lambda\neq 0, \lambda\neq 1$. In this case the genus of the curve is 1 \cite[Theorem 7.2]{Fl00}.

The automorphism $\id\times\iota$ induces the non-symplectic involution $(y:t:x:s)\mapsto (-y:t:x:s)$ on the surface $W$ whose fixed locus consists of the curve $x^{11}s-t(t-s^4)(t-\lambda s^4)\subset W\mathbb{P}(4,1,1)$ whose genus is 10, if $\lambda\neq 0$, $\lambda\neq 1$ \cite[Theorem 7.2]{Fl00}.


\section{Moduli of K3 surfaces}\label{sec: VHS}

By Theorems \ref{theo: PQ with Zp} and \ref{theo: PQ with Z2p} certain components $\mathcal{M}^p_{(n,g(C),k+1)}$ of the moduli space of the K3 surfaces with an automorphism of order $p$ coincide with certain moduli spaces $\mathcal{F}_{PQ}(C_1, D_p)$ of the K3 surfaces which are minimal models either of the quotients $(C_1\times D_p)/({\bf g_1}\times \delta_p)$ or of the quotients $(C_1\times D_p)/({\bf g_1}\times \tau_p)$. Since both $D_p\ra D_p/\delta_p\simeq\mathbb{P}^1$ and $D_p\ra D_p/\tau_p\simeq\mathbb{P}^1$ are branched in 3 points, the parameters of the family depend only on the parameters of the curve $C_1$. In particular the dimension of the family of K3 surfaces is $r-3$, where $r=\sum a_i$ is the number of ramification points of the cover $C_1\ra C_1/{\bf g_1}\simeq \mathbb{P}^1$, cf. Proposition \ref{prop: moduli=r1+r2-6}. Here we describe the relation between the moduli of the curve $C_1$ and the moduli of the surface $S$. In particular we relate the variation of the Hodge structure of weight 2 of  $S$ with the one of $H^1(C_1,\Q)$.

A particular case is the one with $p=3$ and the quotient $(C_1\times D_3)/({\bf g_1}\times \delta_3)$. In this case the variation of the Hodge structure of $S$ (and of a Calabi--Yau 3-fold constructed from $S$ and $D_3$) is described in \cite{GvG}. Moreover, if the family is 1-dimensional, the Picard--Fuchs equation of the surface $S$ is found from the one of the curve $C_1$ (cf. \cite[Section 2.5]{GvG}).

We now assume $S$ to be the minimal model of $(C_1\times D_p)/({\bf g_1} \times \delta_p)$. By construction $S$ admits a non-symplectic automorphism ${\bf g}_S$ induced by $\id\times \delta_p$. The moduli space of the K3 surfaces $S$ obtained in such a way is determined by the variation of the period of $S$ in a certain eigenspace $H^2(S,\C)_{\zeta_p^j}$, cf. Section \ref{sec: non sympl autom}. Indeed, the choice of the period of $S$ determines the Hodge structure of $H^2(S,\C)$ completely.

\begin{prop}\label{prop: half twist} Let $S$ be a generic K3 surface in the family $\mathcal{F}_{PQ}(C_1\times D_p)$ of the surfaces minimal models of $(C_1\times D_p)/({\bf g_1}\times \delta_p)$ listed in \verb|Table 1|. The weight 2 Hodge structure of the transcendental lattice of $S$, $T_S\otimes \Q$, is induced by the weight 1 Hodge structure of $H^1(C_1,\Q)$. In particular, the half twist $(T_S\otimes \Q)_{1/2}$ is $H^1(C_1,\Q)$. As a consequence the dimension of the family $\mathcal{F}_{PQ}(C_1\times D_p)$ is $2g(C_1)/(p-1)-1$.\end{prop}
\proof Since $S$ is generic, the transcendental lattice of $S$ carries a weight 2 Hodge structure of type $(1,(p-1)(m+1)-2,1)$. Since the K3 surface $S$ admits a non-symplectic automorphism ${\bf g}_S$, the Hodge structure of $T_S\otimes\Q$ is of CM-type with the field $K\simeq \Q(\zeta_p)$ (cf. \cite{vG}).  In order to perform a half twist on the Hodge structure one has to fix a CM-type, i.e. a set $\Sigma$ of $(p-1)/2$  distinct embeddings of $K$ in $\C$ with the property that no two of them are conjugate. By abuse of notations we put $\Sigma=\{\zeta_p,\ldots ,\zeta_p^{(p-1)/2}\}$.\\ 
The eigenspaces decomposition (for the action of $\delta_p$) of $H^1(D_p)$ consists of $p-1$ 1-dimensional vector spaces. Therefore $K_{-1/2}\simeq H^1(D_p)$ as Hodge structure of weight 1, where $K_{-1/2}$ is the negative half twist of $K$ (see \cite[Section 1.4]{vG}).\\
Let us denote by $\nu:C_1\times D_p\ra S$ the map induced by the quotient map. The pull-back $\nu^*$ maps $T_S\otimes \Q$ in the $({\bf g_1}\times \delta_p)$ -invariants in $H^1(C,\Q)\otimes H^1(D_p,\Q)$. For dimensional reason $$T_S\otimes\Q\simeq \nu^*(T_S\otimes \Q)=\left(H^1(C,\Q)\otimes_{\Q} H^1(D_p,\Q)\right)^{{\bf g_1}\times \delta_p}.$$
Let us consider the half twist of both the members of the above equation:
\begin{equation}\label{eq: half twist 1}(T_S\otimes\Q)_{1/2}\simeq(\left(H^1(C,\Q)\otimes_{\Q} H^1(D_p,\Q)\right)^{{\bf g_1}\times \delta_p})_{1/2}.\end{equation}
In order to compute the second member of \eqref{eq: half twist 1}, we first consider the $({\bf g_1}\times \delta_p)$-invariant subspace of $H^{1}(C_1,\C)\otimes H^1(D_p,\C)\subset H^{2}(C_1\times D_p)$.  We recall that $H^{1,0}(D_p)_{\zeta_p^i}$ is an eigenspace of dimension 1 if $i\leq (p-i)/2$ and is trivial if $i>(p-1)/2$. By the fact that $p_g(S)=1$, the pair $(C_1,{\bf g_1})$ is such that there exists  only one value $\bar{h}$ such that $\bar{h}>(p-1)/2$ and $H^{1,0}(C_1)_{\zeta_p^{\bar{h}}}$ is non zero.  So the Hodge decomposition of an eigenspace $H^1(C_1)_{\zeta_p^k}$ is trivial, in the sense that $H^1(C_1)_{\zeta_p^k}=H^{a,b}(C_1)_{\zeta_p^k}$ with $(a,b)\in\{(1,0),(0,1)\}$, except for $k=\bar{h}$.  
Hence we obtain:
$$
\begin{array}{ll}&(\left(H^1(C,\C)\otimes H^1(D_p,\C)\right)^{{\bf g_1}\times \delta_p})=\sum_{i=1}^{p-1}\left((H^1(C_1,\C)_{\zeta_p^i}\otimes H^1(D_p,\C)_{\zeta_p^{p-i}}\right)=\\ =& H^{1,0}(C_1)_{\zeta_p^{\bar{h}}}\otimes H^{1,0}(D_p)_{\zeta_p^{p-\bar{h}}}\oplus \sum_{i=1}^{(p-1)/2}\left( H^{1,0}(C_1)_{\zeta_p^{i}}\otimes H^{0,1}(D_p)_{\zeta_p^{p-i}}\right)\oplus\\
\oplus &H^{0,1}(C_1)_{\zeta_p^{p-\bar{h}}}\otimes H^{0,1}(D_p)_{\zeta_p^{\bar{h}}}\oplus \sum_{i=(p+1)/2}^{(p-1)}\left( H^{0,1}(C_1)_{\zeta_p^{i}}\otimes H^{1,0}(D_p)_{\zeta_p^{p-i}}\right).\end{array}
$$
Now we consider the splitting given by the choice of $\Sigma=\{\zeta_p,\ldots \zeta_p^{(p-1)/2}\}$    
and we recall that the action of ${\bf g}_S$ on $T_S$ is induced by the action of $\delta_p$ on $D_p$.
So 
$\left(\sum_{i=(p-1)/2}^{p-1}\left(H^{0,1}(C_1)_{\zeta_p^i}\otimes H^{1,0}(D_p)_{\zeta_p^{p-i}}\right)\right)_{\bar{\Sigma}}=0.$
Hence
$$(\left(H^1(C,\Q)\otimes_{\Q} H^1(D_p,\Q)\right)^{{\bf g_1}\times \delta_p})_{1/2}^{1,0}=\left(
H^{1,0}(C_1)_{\zeta_p^{h}}\otimes H^{1,0}(D_p)_{\zeta_p^{p-h}}\right)\oplus \left(\sum_{i=1}^{(p-1)/2} H^{1,0}(C_1)_{\zeta_p^{i}}\otimes H^{0,1}(D_p)_{\zeta_p^{p-i}}\right).$$
Since $H^{1,0}(D_p)_{\zeta_p^{i}}\simeq \C$ if $i\leq (p-1)/2$, we obtain
$$\left(\left(H^1(C,\Q)\otimes_{\Q} H^1(D_p,\Q)\right)^{{\bf g_1}\times \delta_p}\right)_{1/2}^{1,0}\simeq H^{1,0}(C_1)_{\zeta_p^{h}}\oplus \sum_{i=1}^{(p-1)/2} H^{1,0}(C_1)_{\zeta_p^{i}}=H^{1,0}(C_1).$$
By conjugacy, $\left(\left(H^1(C_1,\Q)\otimes_{\Q} H^1(D_p,\Q)\right)^{{\bf g_1}\times \delta_p}\right)_{1/2}^{0,1}=H^{0,1}(C_1)$. Substituting in \eqref{eq: half twist 1}, we obtain $(T_S\otimes\Q)_{1/2}\simeq H^1(C_1,\C)$ as Hodge structure. \\
In particular $\rk(T_S\otimes\Q)=2g(C_1)$ and by the computation of the moduli of the K3 surface $S$ with a non-symplectic automorphism of order $p$, it follows $m=2g(C_1)/(p-1)-1$.

\endproof

\begin{rem}
{\rm The eigenspaces decomposition (for the action of ${\bf g}_1$) of $H^1(C_1,\Q)$ splits this space in $(p-1)$ equidimensional subspaces and thus $(p-1)(\dim(H^{1}(C_1)_{\zeta_p^j}))=2g(C_1)$ for every $j\in\{1,\ldots, p-1\}$. By the previous proposition we have $\dim(H^{1}(C_1)_{\zeta_p^j})=2g(C_1)/(p-1)=m+1$ for every $j\in\{1,\ldots, p-1\}$.  Moreover, in the proof of the previous proposition we saw that there exists a unique value $\bar{h}$ such that $\bar{h}>p-1$ and $H^{1,0}(C_1)_{\zeta_p^{\bar{h}}}$ is non zero. In particular $\dim(H^{1,0}(C_1)_{\zeta_p^{\bar{h}}})=1$, because $p_g(S)=1$. So $m=\dim(H^1(C_1)_{\zeta_p^{\bar{h}}})-1=\dim(H^{1,0}(C_1)_{\zeta_p^{\bar{h}}}\oplus H^{0,1}(C_1)_{\zeta_p^{\bar{h}}})-1=\dim(H^{1,0}(C_1)_{\zeta_p^{p-\bar{h}}})=\alpha_{p-\bar{h}}$.}
\end{rem}

\begin{rem}{\rm In case $S$ is the minimal model of $(C_1\times D_p)/({\bf g_1}\times \tau_p)$ one can obtain a result similar to the one of Proposition \ref{prop: half twist}: the half twist of $T_S\otimes \Q$ is a sub-Hodge structure of $H^1(C_1,\Q)$ and in fact the one of $\sum_{i=1}^{p-1}H^1(C_1,\C)_{-\zeta_p^i}$. An explicit example is given in \cite[Section 3.11]{vG}.}\end{rem}

The variation of the period of $S$ is described by the Picard--Fuchs equation of $\omega_S$, and so by the Picard--Fuchs equation of certain holomorphic 1-form on $C_1$. In particular if $C_1$ varies in a 1-dimensional family, then it admits an equation of type $y^N=t^a(t-1)^b(t-\lambda)^b$. The forms of these curves and their Picard--Fuchs equations are described in \cite[Section 2.5]{GvG} and this immediately gives the Picard--Fuchs equations of $S$.

For example the Picard-Fuchs equation of the 1-dimensional family $\mathcal{M}^5_{(10,0,2)}$ is the Picard--Fuchs equation of the 1-holomorphic form $\omega_C$ of $C_1$ such that ${\bf g_1}(\omega_C)=\zeta_p^3 \omega_C$.
Since an equation for $C_1$ is $y^5=t(t-1)(t-\lambda)$, the holomorphic form we are interested in is $\omega_1:=dt/y^2$
and its Picard Fuchs equation is
$$\lambda(1-\lambda)\frac{\partial^2 }{\partial\lambda^2}+(\frac{4}{5}-\frac{8}{5}\lambda)\frac{\partial}{\partial\lambda}-\frac{2}{25}=\frac{2}{5}\frac{t(t-1)}{(t-\lambda)y^2}.$$


\bigskip

\bigskip

\textbf{Alice Garbagnati}, \textbf{Matteo Penegini}

\medskip

Dipartimento di Matematica \emph{``Federigo Enriques''}, Universit\`{a} degli Studi di Milano, Via Saldini 50, 20133 Milano, Italy \\

\emph{E-mail addresses:}\\
\verb|alice.garbagnati@unimi.it| \\
 \verb|matteo.penegini@unimi.it|  \\

\newpage
\section*{Appendix. The MAGMA script}
 {\tiny
 \begin{minipage}[t]{8cm}
\begin{verbatim}
Q:=Rationals();
Z:=Integers();

///////////////////////////////////////////

MaxGenus:=function(G)

// Given a cyclic group of order n =  p or
// 2*p it returns the max genus
//  of a curve C with group of
// automorphisms Z/n, C/Z/n=P1, and
// such that the induced action on H^{1,0}(C)
// has an eigenspace of dimension 1.

// Check if the data are correct.

 if not IsCyclic(G) then
  return 0;
 end if;
 if IsCoercible(Z, #G/2) then
  x:=#G/2;
  y:=Z ! x;
  if not IsPrime(y) then
  return 0;
 end if; end if;
 if not IsCoercible(Z, #G/2) then
  if not IsPrime(#G) then
  return 0;
 end if; end if;
 R:=PolynomialRing(Q,#G-1);

 Preparation:=function(G)

// Given a cyclic group G of order
// 2*p or p it return the matrix of hol.
// Lefschetz.
// WARNING: the eigenspaces have a peculiar
// ordering. The integers are always first.

  T:=CharacterTable(G);
  if (#G mod 2) eq 0 then
   x:=#G/2;
   y:=Z ! x;
   F:=CyclotomicField(y);
  else
   y:=#G;
   F:=CyclotomicField(#G);
  end if;
  L:=[];
  Append(~L,T[#G]);
if #G ge 5 then
  Append(~L,T[#G-3]);
  Append(~L,T[2]);
else
  Append(~L,T[2]);
  Append(~L,T[1]);
end if;
  R:=PolynomialRing(F,#G-1);
  List:=[];
  for i in [1..#G-1] do
   List[i]:=R.i;
  end for;
  v:=Vector(R,List);
  g:=[];
  for j in [1..3] do
   List2:=Eltseq(L[j]);
   w:=Vector(R,List2);
   g[j]:=&+[w[i+1]*v[i]: i in [1..#G-1]];
  end for;
  Au:=[];
  for k in [1..3] do
   C:=Coefficients(g[k]);
   TN:=[];
   for hh in [1..y-1] do
    TN[hh]:=[];
   end for;
   for i in [1..#G-1] do
    pip:=Eltseq(C[i]);
    for h in [1..y-1] do
     Append(~TN[h],-pip[h]);
    end for;
   end for;
   for kk in [1..y-1] do
    Append(~Au,TN[kk]);
   end for; end for;
\end{verbatim}
\end{minipage}
\begin{minipage}[t]{8cm}
\begin{verbatim}
  Li:=[];
  for i in [1..#G-1] do
   for j in [1..#G-1] do
    Append(~Li,Au[i][j]);
  end for; end for;
  M:=Matrix(R,#G-1,Li);
  return Li;
 end function;

 TerminiNoti:=function(G)
  R:=PolynomialRing(Q,#G-1);
  if IsPrime(#G) then
   p:=#G;
  else
   x:=(#G)/2;
   p:=Z ! x;
  end if;
  sf:=(p-1)/2;
  st:=Z ! sf;
  F<s>:=CyclotomicField(p);

// First case |G|=p

  if IsPrime(#G) then
   Li:=[];
   for i in [1..p-1] do
    Li[i]:=1/(1-s^i);
   end for;
   P:=[];
   List:=[];
   for i in [1..p-1] do
    List[i]:=Eltseq(Li[i]);
   end for;
   K<s>:=PolynomialRing(Q,1);
   for i in [1..p-1] do
    P[i]:=&+[List[i][j]*s^(j-1): j in [1..p-1]];
   end for;
   ww:=[i*0: i in [1..p-1]];
   w:=Matrix(R,1,p-1,ww);
   for i in [1..p-1] do
    w:=w+R.i*Matrix(R,1,p-1,Coefficients(P[i]));
   end for;
   v:=[];
   v[1]:=-1+w[1,p-1];
   for i in [2..p-1] do
    v[i]:=w[1,p-i];
   end for;
   return v;
  else

// Second case |G|=2*p

   Li:=[];
   for i in [1..p-1] do
    Li[i]:=1/(1+s^i);
   end for;
   P:=[];
   List:=[];
   for i in [1..p-1] do
    List[i]:=Eltseq(Li[i]);
   end for;
   TN:=[];
   for hh in [1..p-1] do
    TN[hh]:=[];
   end for;
   for i in [1..p-1] do
    for k in [1..p-1] do
     Append(~TN[k],List[i][k]);
   end for;end for;
   w:=[];
   for i in [1..p-1] do
    w[i]:=&+[TN[i][j]*R.j: j in [1..p-1]];
   end for;
   Li:=[];
   for i in [1..p-1] do
    Li[i]:=1/(1-s^i);
   end for;
   P:=[];
   List:=[];
   for i in [1..p-1] do
    List[i]:=Eltseq(Li[i]);
   end for;
   TN:=[];
    \end{verbatim}
\end{minipage}
\newpage
\begin{minipage}[t]{8cm}
\begin{verbatim}
   for hh in [1..p-1] do
    TN[hh]:=[];
   end for;
   for i in [1..p-1] do
    for k in [1..p-1] do
     Append(~TN[k],List[i][k]);
   end for;end for;
   ww:=[];
   for i in [1..p-1] do
    ww[i]:=(&+[TN[i][j]*R.((j+j*st) mod p): j in [1..p-1]])
+(&+[TN[i][j]*R.(((j+j*st) mod p)+p-1): j in [1..p-1]]);
   end for;
   t:=2*p-1;
   www:=(&+[R.i: i in [1..p-1]]+R.t)/2;
   v:=[];
   for i in [1..p-1] do
    if i eq 1 then
     v[i]:=-1+w[i];
    else
     v[i]:=w[i];
    end if;
   end for;
   for i in [1..p-1] do
    if i eq 1 then
     v[p-1+i]:=-1+ww[i];
    else
     v[p-1+i]:=ww[i];
   end if; end for;
   v[2*p-1]:=www-1;
  end if;
  return v;
 end function;

// Find the correct eigenspace, which will
// have dimension 1

 if IsPrime(#G) then
  if #G ge 7 then
   p1:=1;
  else p1:=#G-1;
  end if;
 end if;
 if not IsPrime(#G) then
  if #G ge 7 then
   x:=(#G)/2+1;
   p1:=Z ! x;
  else
   p1:=#G-1;
  end if;
 end if;

// MAIN RUTINE MAXGENUS

 z:=#G-1;
 Li:=Preparation(G);
 v:=TerminiNoti(G);
 M:=Matrix(R,#G-1,Li);
 v1:=Vector(R,v);
 M1:=Transpose(M);
 sol:=Solution(M1,v1);
 g1:=&+[sol[i]: i in [1..z]];
 printf "========================== \n";
 printf "The general genus is: \n";
 printf "%o, \n", g1;
 printf "========================== \n";
 printf "\n";
 printf "========================== \n";
 printf "The dim of eigsp. are:  \n";
 for i in [1..z] do
 printf" %o, \n", sol[i];
 end for;
 printf "========================== \n";
 Coeff:=Coefficients(sol[p1]-1);
 ass:=-1/Coeff[1]*(sol[p1]-1-Coeff[1]*R.1);

// Get a bound on g
 g:=Evaluate(g1,1,ass);
 printf "\n";
 printf "========================== \n";
 printf "The special genus is: \n";
 printf "%o, \n", g;
 printf "========================== \n";
 printf "\n";
 Coefs:=Coefficients(g);
 MaxGen:=Coefs[#G-1];
 \end{verbatim}
\end{minipage}
\begin{minipage}[t]{8cm}
\begin{verbatim}
 printf "========================== \n";
 printf "The max genus is: %o \n", MaxGen;
 printf "========================== \n";
 printf "\n";
 return MaxGen,sol;
end function;

/////////////////////////////////////////////
//
//             END OF MAXGENUS
//
/////////////////////////////////////////////

// THE FOLLOWING PART WORKS ONLY FOR THE GROUP Z/pZ

Nram:=function(MaxGen,p)

// Given the genus g(C) genus and a prime p
// returns the number of branch points
// of the cover C -> C/Z/p = P1

 x:=2*MaxGen-2;
 M:=[p];
 y:=p*(-2+#M-1/M[1]);
 while x-y gt 0  do
  Append(~M, p);
  z:=#M/p;
  y:=p*(-2+#M-z);
 end while;
return #M;
end function;

PossRami:=function(N,NN)

// It returns a seq of seqs of
// #NN, whos sum of elements
// is less or eq to N

 S:={1..N};
 M:=[];
 i:=1;
  while i le N do
   for h in [1..NN] do
    Append(~M,RestrictedPartitions(i,h,S));
   end for;
   i:=i+1;
  end while;
  K:=[];
  for LL in M do
   for i in [1..#LL] do
    L:=[];
    for j in [1..#LL[i]] do
     L[j]:=LL[i][j];
    end for;
    for k in [#LL[i]+1..NN] do
     L[k]:=0;
    end for;
    Append(~K,L);
   end for;
  end for;
  return K;
end function;

MaybeSur1:=function(s,p)

// Given vectors in H return all possible pairings
// of dim of eigsp. which could give a surface with
// pg=1

 PP:=[];
// all possible pairing
 for i in [1..#s] do
  for j in [i..#s] do
   Append(~PP,[s[i],s[j]]);
 end for;end for;
 GP:=[];
 for i in [1..#PP] do
  n:=0;
  for j in [1..2] do
   for k in [1..p-1] do
    if PP[i][j]`EignSpaces[k] eq 0 then n:=n+1;
  end if;end for;end for;
  if n ge p-2 then Append(~GP, PP[i]);
 end if;end for;
 return GP;
end function;
 \end{verbatim}
\end{minipage}
\newpage
\begin{minipage}[t]{8cm}
\begin{verbatim}
TheSur:=function(MaybeSur,p,RF,sol)

 Sur:=[];
 for j in [1..#MaybeSur] do
  for h in [1..p-1] do
   v1:=Rotate(MaybeSur[j][2]`EignSpaces,h);
   if &+[MaybeSur[j][1]`EignSpaces[k]*v1[k]: k in [1..p-1]]
   eq 1 then
    for t in [1..p] do
     v2:=Rotate(MaybeSur[j][2]`FixPoints,t);
     if IsCoercible(Z,Evaluate(sol[1],v2)) then
      w:=rec< RF | EignSpaces := v1, FixPoints:=v2>;
      Include(~Sur,[MaybeSur[j][1],w]);
 end if;end for;end if;end for;end for;
 return Sur;
end function;

///////////////////////////////////////////////////
//
// Spherical Systems of Generators
// for groups of order p or 2p
//
///////////////////////////////////////////////////

GenSys:=function(G,m);

// Given a group G of order n it returns
// m integers whose sum is a multiple of n.

 n:=#G;
 H:=[];
 L:=[];
 M:={};
 for i in [1..n-1] do
  Include(~M,i);
 end for;
 for x in [1..m-1] do
  y:=n*x-(n*x mod n);
  TL:=RestrictedPartitions(y,m,M);
  for j in [1..#TL] do
   Append(~L,Reverse(TL[j]));
  end for;
 end for;
 return L;
end function;

VecRami:=function(H,G,n)

// Given a cyclic group of order p or 2p it returns
// a list of spherical system of generators of size n

 K:=[];
 if IsPrime(#G) then
  g:=G.1;
 else
  g:=G.1*G.2;
 end if;
  for h in H do
   L:=[];
   for i in [1..n] do
    Append(~L,g^h[i]);
   end for;
   if #sub<G | L> eq #G then
    Append(~K, L);
   end if;
  end for;
 return K;
end function;

///////////////////////////////////////////////////
//
// Match the Eigenspaces such that p_g of the surface is 1
//
///////////////////////////////////////////////////

GioCop:=function(H,K,G);
 sum:=0;
 p:=#G;
 for i in [1..p-1] do
  sum:=sum+H[i]*K[p-i];
 end for;
 return sum;
end function;
 \end{verbatim}
\end{minipage}
\begin{minipage}[t]{8cm}
\begin{verbatim}
// Reordering the Eigenspaces only for group
// of order 2p
//
///////////////////////////////////////////////////
EigReorderingp:=function(V,G,sol)

// The eigenspaces are reordered
// w.r.t. -z, -z^2, -z^3, ...

 p:=Z ! #G;
 L:=[];
 for i in [1..p-1] do
  L[i]:=p*Coefficients(sol[i])[1];
 end for;
 NewV:=[i*0: i in [1..#G-1]];
 for i in [1..p-1] do
  NewV[Z ! L[i]]:=NewV[Z ! L[i]]+V[i];
 end for;
 return NewV;
end function;


EigReordering:=function(V,G,sol)

// The eigenspaces are reordered
// w.r.t. -z, z^2, -z^3, ...

 x:=(#G)/2;
 p:=Z ! x;
 q1:=(p-1)/2;
 q:= Z ! q1;
 L:=[];
 for i in [1..p-1] do
  L[i]:=2*p*Coefficients(sol[p+i])[1];
 end for;
 NewV:=[i*0: i in [1..#G-1]];
// EigS rel to elm of order 2
 NewV[p]:=NewV[1]+V[1];
// EigSs rel to elms of order p
 for i in [1..p-1] do
  if (Z ! L[i]) le p then
   NewV[(Z !L[i])+1]:=NewV[(Z ! L[i])+1]+V[i+1];
  else
   NewV[(Z ! L[i])-1]:=NewV[(Z ! L[i])-1]+V[i+1];
  end if;
 end for;
// EigSs rel to elms of order 2p
 for i in [1..p-1] do
  NewV[Z ! L[i]]:=NewV[Z ! L[i]]+V[p+i];
 end for;
 return NewV;
end function;

///////////////////////////////////////////////////
//
// Fix points reordering
//
///////////////////////////////////////////////////

FixPointsReorderingp:=function(V, G, tmp)

 F:=FiniteField(#G);
 g:=G.1;
 Gelm:=[];
 for i in [1..#G] do
  Append(~Gelm, g^i);
 end for;
 FixPoints:=[i*0: i in [1..#G-1]];
 for a in [1..tmp] do
  x:=F ! Position(Gelm, G ! V[a]);
  y1:=x^-1;
  y:= Z ! y1;
  FixPoints[y]:=(FixPoints[y]+1);
 end for;
 return FixPoints;
end function;


FixPointsReordering:=function(V, G, tmp)

// Needed for the calculation of the dim
// of the eigenspaces

 p1:=#G/2;
 p:= Z !p1;
  \end{verbatim}
\end{minipage}
\newpage
\begin{minipage}[t]{8cm}
\begin{verbatim}
 F:=FiniteField(p);
 g:=G.1*G.2;
 Gelm:=[];
// Fixing the one generator G.1*G.2=z we list z,z^2,z^3..
 for i in [1..#G] do
  Append(~Gelm, g^i);
 end for;
 OrdGelm:=[];
 for h in Gelm do
  Append(~OrdGelm, Order(h));
 end for;
 FixPoints:=[i*0: i in [1..#G-1]];
 for a in [1..tmp] do
  if Position(Gelm,V[a]) eq p then
   FixPoints[#G-1]:=FixPoints[#G-1]+p;
  end if;
   if IsEven(Position(Gelm,V[a])) then
    x1:=(Position(Gelm, V[a]))/2;
    x:= F ! x1;
    x2:=x^-1;
    x3:= Z ! x2;
    y:= Z ! p-1+x3;
    FixPoints[y]:=FixPoints[y]+2;
   end if;
    if IsOdd(Position(Gelm,V[a])) then
     if Position(Gelm,V[a]) ne p then
      x:=F ! Position(Gelm, G ! V[a]);
      y1:=x^-1;
      y:= Z ! y1;
      FixPoints[y]:=(FixPoints[y]+1);
     end if;
    end if;
 end for;
 return FixPoints;
end function;

FixPointsReordering2:=function(V, G, tmp)

// The fixpoints are ordered as [-z,-z^2,..][z,z^2,..][z^p]
// V is an ssg, tmp is the number of ramification pts.

 p1:=#G/2;
 p:= Z !p1;
 F:=FiniteField(p);
 g:=G.1*G.2;
 Gelm:=[];
 for i in [1..#G] do
  Append(~Gelm, g^i);
 end for;
 OrdGelm:=[];
 for h in Gelm do
  Append(~OrdGelm, Order(h));
 end for;
 FixPoints:=[i*0: i in [1..#G-1]];
 for a in [1..tmp] do
  if Position(Gelm,V[a]) eq p then
   FixPoints[#G-1]:=FixPoints[#G-1]+p;
  end if;
  if IsEven(Position(Gelm,V[a])) then
   x1:=(Position(Gelm, V[a]))/2;
   x:= F ! x1;
   x2:=x^-1;
   x3:=x2*2;
   x4:= Z ! x3;
   y:= Z ! p-1+x4;
   FixPoints[y]:=FixPoints[y]+2;
  end if;
  if IsOdd(Position(Gelm,V[a])) then
   if Position(Gelm,V[a]) ne p then
    x:=F ! Position(Gelm, G ! V[a]);
    y1:=x^-1;
    y:= Z ! y1;
    FixPoints[y]:=(FixPoints[y]+1);
   end if;
  end if;
 end for;
 return FixPoints;
end function;

///////////////////////////////////////////////////
//
// We couple the curves C_1 and C_2 to return
// potential surfaces of the form T=(C_1 x C_2)/Z/(2p)Z
// such that p_g(T)=1.
//
///////////////////////////////////////////////////
  \end{verbatim}
\end{minipage}
\begin{minipage}[t]{8cm}
\begin{verbatim}
Surface:=function(GT,GT3,sol,G,t1,t2)

 RF := recformat<SSG, EignSpaces : SeqEnum, FixPoints : SeqEnum >;
 Sur:=[];
 if IsPrime(#G) then
  for i in [1..#GT] do
   W21b:=[];
   for k in [1..#G-1] do
    C:=FixPointsReorderingp(GT[i],G,t1);
    W21b[k]:=Evaluate(sol[k],C);
   end for;
   W21:=EigReorderingp(W21b,G, sol);
   for j in [i..#GT3] do
    W22b:=[];
    for k in [1..#G-1] do
     C:=FixPointsReorderingp(GT3[j],G,t2);
     W22b[k]:=Evaluate(sol[k],C);
    end for;
    W22:=EigReorderingp(W22b,G, sol);
    if GioCop(W21,W22,G) eq 1 then
    s:=[];
    Append(~s, rec< RF | SSG:= GT[i], EignSpaces := W21, FixPoints:=
     FixPointsReorderingp(GT[i],G,t1)>);
    Append(~s, rec< RF | SSG:= GT3[j], EignSpaces := W22, FixPoints:=
     FixPointsReorderingp(GT3[j],G,t2)>);
    Append(~Sur, s);
   end if;
  end for;
 end for;
 else
  for i in [1..#GT] do
   W21b:=[];
   for k in [1..#G-1] do
    C:=FixPointsReordering(GT[i],G,t1);
    W21b[k]:=Evaluate(sol[k],C);
   end for;
   W21:=EigReordering(W21b,G, sol);
   for j in [i..#GT3] do
    W22b:=[];
    for k in [1..#G-1] do
     C:=FixPointsReordering(GT3[j],G,t2);
     W22b[k]:=Evaluate(sol[k],C);
    end for;
    W22:=EigReordering(W22b,G, sol);
   if GioCop(W21,W22,G) eq 1 then
    s:=[];
    Append(~s, rec< RF | SSG:= GT[i], EignSpaces := W21, FixPoints:=
    FixPointsReordering2(GT[i],G,t1)>);
    Append(~s, rec< RF | SSG:= GT3[j], EignSpaces := W22, FixPoints:=
    FixPointsReordering2(GT3[j],G,t2)>);
    Append(~Sur, s);
   end if;
  end for;
 end for;
 end if;
 return Sur;
end function;

///////////////////////////////////////////////////
//
// Singularities check is needed. First we calculate the contributions
// of each singularity to k^2_S, e(S), and \chi(S). We borrowed this part
// of the program from [BP11].
//
///////////////////////////////////////////////////

ContFrac:=function(s)
  CF:=[ ]; r:=1/s;
  while not IsIntegral(r) do
    Append(~CF, Ceiling(r)); r:=1/(Ceiling(r)-r);
  end while;
  return Append(CF, r);
end function;

Nq:=func<cf|#cf eq 1 select cf[1] else cf[1]-1/$$(Remove(cf,1))>;

RatNum:=func<seq|1/Nq(seq)>;

// "Wgt" computes the weight of a sequence, i.e., the sum of its
// entries. It bounds strictly from below B of the corresponding
// singular point.

Wgt:=function(seq)
  w:=0; for i in seq do w+:=i; end for; return w;
end function;
  \end{verbatim}
\end{minipage}
\newpage
\begin{minipage}[t]{8cm}
\begin{verbatim}
// The next script computes all rational number
// whose continuous fraction has small weight,
// by listing all sequences (modulo
// "reverse") and storing the corresponding rational number.

RatNumsWithSmallWgt:=function(maxW)
  S:={ }; T:={}; setnums:={RationalField()| };
  for i in [2..maxW] do Include(~S, [i]); end for;
  for i in [1..Floor(maxW/2)-1] do
  for seq in S do
    if #seq eq i then
    if maxW-Wgt(seq) ge 2 then
    for k in [2..maxW-Wgt(seq)] do
     Include(~S,Append(seq, k));
    end for; end if; end if;
  end for; end for;
  for seq in S do
  if Reverse(seq) notin T then Include(~T,seq);
  end if; end for;
  for seq in T do Include(~setnums, RatNum(seq)); end for;
  return setnums;
end function;

// The next two scripts compute the invariants
// B and e of a rational number (i.e., of
// the corresponding singular point).

InvB:=func<r|Wgt(ContFrac(r))+r+RatNum(Reverse(ContFrac(r)))>;

Inve:=func<r|#ContFrac(r)+1-1/Denominator(RationalField()!r)>;

// Here is the invariant k of the basket:

Invk:=func<r|InvB(r)-2*Inve(r)>;

kappa:=function(p)
 KK:=[];
 for i in [1..p-1] do
  KK[i]:=-Invk(i/p);
 end for;
 return KK;
end function;

EuleR:=function(p)
 CC:=[];
 for i in [1..p-1] do
  CC[i]:=Inve(i/p);
 end for;
 return CC;
end function;

///////////////////////////////////////////////////
//
// The program for singularities for surfaces with
// group Z/pZ
//
///////////////////////////////////////////////////

TypeSing:=function(Sur, p)

 count:=0;
 F:=FiniteField(p);
 KK:=kappa(p);
 CC:=EuleR(p);
 for k in [1..#Sur] do
  L:=[i*0: i in [1..p-1]];
  for i in [1..(p-1)] do
   for j in [1..(p-1)] do
    x:= F ! i;
    y:= F ! j;
    z1:=x*y^-1;
    z:= Z ! z1;
    if z ne 0 then
     L[z]:=L[z]+Sur[k][1]`FixPoints[i]*Sur[k][2]`FixPoints[j];
     g1:=&+[Sur[k][1]`EignSpaces[h]: h in [1..p-1]];
     g2:=&+[Sur[k][2]`EignSpaces[h]: h in [1..p-1]];
    end if;
   end for;
  end for;
  K2:=8*(g1-1)*(g2-1)/p+&+[L[t]*KK[t]: t in [1..p-1]];
  chi:=4*(g1-1)*(g2-1)/p+&+[L[t]*CC[t]: t in [1..p-1]]+K2;
  if chi eq 24 then
   printf "========================== \n";
   printf"there is a surface with curves %o \n", Sur[k];
   printf"with singularities: \n" ;
   for h in [1..p-1] do
    printf"%o x 1/%o(1,%o) \n", L[h],p,h;
      \end{verbatim}
\end{minipage}
\begin{minipage}[t]{8cm}
\begin{verbatim}
   end for;
   printf "the minimal resolution has at least %o -1-curves \n", -K2;
   printf "the minimal model has 12*chi = %o \n", chi;
   printf "========================== \n";
   printf "\n";
   count:=count+1;
  end if;
 end for;
 return count;
end function;

///////////////////////////////////////////////////
//
// The program for singularities for surfaces with
// group Z/2pZ
//
///////////////////////////////////////////////////

Singulp:=function(V,W,p)
 F:=FiniteField(p);
 L:=[i*0: i in [1..p-1]];
 for i in [1..(p-1)] do
  for j in [1..(p-1)] do
   x:= F ! i;
   y:= F ! j;
   z1:=x*y^-1;
   z:= Z ! z1;
   if z ne 0 then
    L[z]:=L[z]+V[i]*W[j];
 end if; end for; end for;
 return L;
end function;

Singul2p:=function(V,W,p)
 F:=FiniteField(p);
 L:=[i*0: i in [1..p-1]];
 for i in [1..(p-1)] do
  for j in [1..(p-1)] do
   x:= F ! 2*i;
   y:= F ! j;
   z1:=x*y^-1;
   z:= Z ! z1;
   if z ne 0 then
    L[z]:=L[z]+V[i]*W[j];
 end if; end for; end for;
 return L;
end function;

Singul2p2:=function(V,W,p)
 F:=FiniteField(p);
 L:=[i*0: i in [1..p-1]];
 for i in [1..(p-1)] do
  for j in [1..(p-1)] do
   x:= F ! 2*i;
   y:= F ! j;
   z1:=x^-1*y;
   z:= Z ! z1;
   if z ne 0 then
    L[z]:=L[z]+V[i]*W[j];
 end if;end for;end for;
 return L;
end function;


///////////////////////////////////////////////////
//
// We check if the surfaces given in Sur are potential
// K3 calculating the invariants of (C_1 x C_2)/G
//
///////////////////////////////////////////////////

CheckSing:=function(Sur,G)
 x:=#G/2;
 p:= Z ! x;
 count:=0;
 Super:=[];
 Sing:=[];
 os:=#Sing;
 Singolar:=[];
 for k in [1..#Sur] do
  VA:=[];
  WA:=[];
  VB:=[];
  WB:=[];
    \end{verbatim}
\end{minipage}
\newpage
\begin{minipage}[t]{8cm}
\begin{verbatim}
  for i in [1..p-1] do
   VA[i]:=Sur[k][1]`FixPoints[i];
   WA[i]:=Sur[k][2]`FixPoints[i];
   VB[i]:=Sur[k][1]`FixPoints[i+p-1];
   WB[i]:=Sur[k][2]`FixPoints[i+p-1];
  end for;

// Calculate the singularities divided in 3 types
// first the ones of the form 1/#G (1,*)
// then the ones of the form 1/p (1,*)
// and in the end of the form 1/2 (1,1)

  Singolar[1]:=Singulp(VA,WA,p);
  Singolar[2]:=[];
  for t in [1..p-1] do
   Append(~Singolar[2], Z ! (Singulp(VB,WB,p)[t]+
   Singul2p(VA,WB,p)[t]+Singul2p2
(WA,VB,p)[t])/2);
  end for;
  Singolar[3]:=
  [Z ! (Sur[k][1]`FixPoints[#G-1]*Sur[k][2]`FixPoints[#G-1]+
&+[Sur[k][2]`FixPoints[#G-1]*Sur[k][1]`FixPoints[h]:
h in [1..p-1]]+
&+[Sur[k][1]`FixPoints[#G-1]*Sur[k][2]`FixPoints[h]:
h in [1..p-1]])/p];

// We calculate the genera of the curves

  g1:=&+[Sur[k][1]`EignSpaces[h]: h in [1..#G-1]];
  g2:=&+[Sur[k][2]`EignSpaces[h]: h in [1..#G-1]];
  Include(~Sing,Singolar);
  ns:=#Sing;

// Calculate K^2

  kcontrn:=[];
  kcontrp:=[];
  RR:=quo< Z | #G>;
  U1:=[];
  for i in [1..#G-1] do
   x:=RR ! i;
   if IsUnit(x) then
    Append(~U1,x);
   end if;
  end for;
  for i in [1..p-2] do
   if IsOdd(i) then
    kcontrn[i]:=Invk((i)/#G)*Singolar[1][i];
   else
    kcontrn[i]:=Invk((i+p)/#G)*Singolar[1][i];
   end if;
   kcontrp[i]:=Invk(i/p)*Singolar[2][i];
  end for;
  kappa2:=8*(g1-1)*(g2-1)/#G-&+[kcontrn[t]: t in [1..p-2]]-
  &+[kcontrp[t]: t in [1..p-2]];

// Calculate e the Euler number

  econtrn:=[];
  econtrp:=[];
  RR:=quo< Z | #G>;
  U:=[];
  for i in [1..#G-1] do
   x:=RR ! i;
   if IsUnit(x) then
    Append(~U,x);
  end if;end for;
  for i in [1..p-1] do
   if IsOdd(i) then
    econtrn[i]:=Inve(i/#G)*Singolar[1][i];
   else
    econtrn[i]:=Inve((i+p)/#G)*Singolar[1][i];
  end if;end for;
  for i in [1..p-1] do
   econtrp[i]:=Inve(i/p)*Singolar[2][i];
  end for;
  econtr12:=Inve(1/2)*Singolar[3][1];
  euler:=4*(g1-1)*(g2-1)/#G+&+[econtrn[t]: t in [1..p-1]]+&
  +[econtrp[t]: t in [1..p-1]]+econtr12;

// Calculate \chi

  chi:=(euler+kappa2)/12;

// Only potential K3s survive the two tests below

        \end{verbatim}
\end{minipage}
\begin{minipage}[t]{8cm}
\begin{verbatim}
  if chi eq 2 then
   if os ne ns then
    Append(~Super,Sur[k]);
    printf "========================== \n";
    printf"there is a surface with curves %o \n", Sur[k];

    printf"with singularities: \n" ;
    for i in [1..p-1] do
     if IsOdd(i) then
      printf"%o x 1/%o(1,%o) \n", Singolar[1][i],#G,i;
     else
      printf"%o x 1/%o(1,%o) \n", Singolar[1][i],#G,(i+p);
    end if;end for;
    for i in [1..p-1] do
     printf"%o x 1/%o(1,%o) \n", Singolar[2][i], p, i;
    end for;
    printf"%o x 1/2(1,1) \n", Singolar[3][1];
    printf "the minimal resolution has at least %o -1-curves \n", -kappa2;
    printf "the minimal model has chi = %o \n", chi;
    printf "========================== \n";
    printf "\n";
    count:=count+1;
  end if; end if;
  os:=#Sing;
 end for;
 return count, Super;
end function;

///////////////////////////////////////////////////
//
// MAIN ROUTINES and COMMANDS
//
///////////////////////////////////////////////////
//
// For ANY surfaces with group Z/(p)Z
//
///////////////////////////////////////////////////
Surfacesp:=function(G,c)

// Set c=true to perform the calculation of Eigenspaces and max
// genus only.
// Set c=false to get a list of possible K3, in this case the
// algorithm is not optimal and for p >=11 could be extremly slow!

 p:=Z ! #G;
 R:=PolynomialRing(Q,p-1);
 MaxGen,sol:=MaxGenus(G);
 Nrami:=Nram(MaxGen,p);
 printf "========================== \n";
 printf "The max n. of rami is: %o \n", Nrami;
 printf "========================== \n";
 printf "\n";
  if c then return 0;
  end if;
 K:=PossRami(Nrami,p-1);
 printf "========================== \n";
 printf "Starting comb calc fom: %o \n", #K;
 printf "It will be very very slow... \n";
 printf "========================== \n";
 printf "\n";
 S:={1..p-1};
 To:=[];
 t:=#K;
 StS:=SetToSequence(Permutations(S));
 for i in [1..t] do
  for j in [1..#StS] do
   C:=[];
   for k in [1..p-1] do
    C[k]:=K[i][StS[j][k]];
   end for;
   Append(~To,C);
 end for;end for;
 Tot:=SetToSequence(Seqset(To));
 printf "========================== \n";
 printf "Finish transposition \n";
 printf "Starting cyclic permut.... \n";
 printf "========================== \n";
 printf "\n";
 printf "========================== \n";
 printf "The n. of Tot. is: %o \n", #Tot;
 printf "========================== \n";
 printf "\n";
 ops:=0;
      \end{verbatim}
\end{minipage}
\newpage
\begin{verbatim}
 RF := recformat< EignSpaces : SeqEnum, FixPoints : SeqEnum >;
 s:=[];
 H:=[];
 printf "========================== \n";
 for i in [1..#Tot] do
  W22:=[];
  for j in [1..p-1] do
   Append(~W22, Evaluate(sol[j],Tot[i]));
  end for;
  OS:=#H;
  if W22[1] eq 1 then
   Include(~H,W22);
   Append(~s, rec< RF | EignSpaces := W22, FixPoints:= Tot[i]>);
   NS:=#H;
   if ops eq 0 then
    if #H eq 1 then
     ops:=ops+1;
    end if;
   end if;
  end if;
 end for;
 printf "========================== \n";
 printf "\n";
 printf "========================== \n";
 printf"We have %o curves \n", #s;
 printf "========================== \n";
 printf "\n";
 MaybeSur:=MaybeSur1(s,p);
 Sur:=TheSur(MaybeSur,p,RF,sol);
 printf "\n";
 printf "========================== \n";
 printf"We have %o possible Surfaces \n", #Sur;
 printf "========================== \n";
 printf "\n";
 c:=TypeSing(Sur, p);
 printf "\n";
 printf "========================== \n";
 printf"We have %o Surfaces \n", c;
 printf "========================== \n";
 printf "\n";
 return Sur;
end function;

///////////////////////////////////////////////////
//
// For t1-t2 points ramifications surfaces
// with group either Z/(p)Z or Z/(2p)Z
//
///////////////////////////////////////////////////

t1t2PtsSurfaces:=function(G,t1,t2)

 R:=PolynomialRing(Q,#G-1);
 MaxGen,sol:=MaxGenus(G);
 GTT:=GenSys(G,t1);
 GT:=VecRami(GTT,G,t1);
 GTT3:=GenSys(G,t2);
 GT3:=VecRami(GTT3,G,t2);
 Sur:=Surface(GT,GT3,sol,G,t1,t2);
 if IsPrime(#G) then
  count:=TypeSing(Sur,#G);
 else
  count, Super:=CheckSing(Sur,G);
 end if;
 printf "\n";
 printf "========================== \n";
 printf"We have %o possible Surfaces \n", count;
 printf "========================== \n";
 printf "\n";
 return count;
end function;
\end{verbatim}
}



\begin{thebibliography} {9}
%
\bibitem[AS08]{AS08}
M. Artebani, A. Sarti \textit{Non-symplectic automorphisms of
order 3 on {$K3$} surfaces}, Math. Ann.,  \textbf{342} (2008), 903--921.
%
\bibitem[AST11]{AST}
M. Artebani, A. Sarti, S Taki, \textit{{$K3$} surfaces with
non-symplectic automorphisms of prime order}, Math. Z.,
\textbf{268}, (2011), 507--533.
%
\bibitem[BHPV]{BHPV}
W.P. Barth, K. Hulek, C.A.M. Peters, A. Van de Ven,
\textit{Compact complex surfaces}. Second edition. Ergebnisse der
Mathematik und ihrer Grenzgebiete. 3. Folge. A Series of Modern
Surveys in Mathematics [Results in Mathematics and Related Areas.
3rd Series. A Series of Modern Surveys in Mathematics],
\textbf{4}. Springer-Verlag, Berlin, 2004.
%
\bibitem[BCGP09]{BCGP}
I. Bauer, F. Catanese, F. Grunewald, R. Pignatelli,
\textit{Quotients of products of curves, new surfaces with $p_g=0$
and their fundamental groups}. American Journal of Mathematics,
\textbf{134}, (2012), 993--1049.
%
\bibitem[BP12]{BP11}
I. Bauer, R. Pignatelli, \textit{The classification of minimal
product-quotient surfaces with $p_g=0$}. Mathematics of
Computation \textbf{81}, (2012), 2389--2418.
%
%
\bibitem[C00]{cat00} 
F. Catanese, \textit{Fibred surfaces, varieties isogenous to a
product and related moduli spaces}. Amer. J. Math. \textbf{122},
(2000), 1--44.
%
\bibitem[CW34]{CW}
C. Chevalley, A. Weil, \textit{\"Uber das Verhalten der Intergrale 1. Gattung bei Automorphismen des Funktionenkorpers}. Abhand. Math. Sem. Hamburg {\bf 10} (1934), 358--361.
%
\bibitem[Di12]{Dillies}
J. Dillies, {\it On some order 6 automorphisms of elliptic K3 surfaces}, Albanian J. Math. {\bf 6} (2012).
%
\bibitem[D82]{D82}
I. Dolgachev, \textit{Weighted projective varieties}. Group actions and vector fields (Vancouver, B.C., 1981), 34--71, Lecture Notes in Math., 956, Springer, Berlin, 1982.
%
\bibitem[DK07]{DK}
I. V. Dolgachev, S. Kond$\bar{\rm o}$, {\it Moduli of K3 surfaces and complex ball quotients},
In Arithmetic and geometry around hypergeometric functions,
{\bf 260} of Progr. Math., 43--100. Birkh \"auser, Basel, 2007.
%
%
\bibitem[Fl00]{Fl00}
A. R. Iano-Fletcher, \textit{Working with weighted complete
intersections}, in Explicit birational geometry of 3-folds, London
Math. Soc. Lecture Note Ser.,    \textbf{281}. Cambridge Univ.
Press. (2000), 101--173.
%
\bibitem[F71] {F71}
E. Freitag, {\it \"{U}ber die {S}truktur der {F}unktionenk\"orper zu
hyperabelschen {G}ruppen. {I}}, J. Reine Angew. Math. {\bf 247} (1971), 97--117.
%
\bibitem[GvG10]{GvG} A. Garbagnati, B. van Geemen, {\it The Picard-Fuchs equation of a family of Calabi-Yau threefolds without maximal unipotent monodromy},  Int. Math. Res. Not. IMRN  {\bf 16}  (2010), 3134--3143.
%
\bibitem[GS13]{GS} A. Garbagnati, A. Sarti {\it Symplectic and non-symplectic automorphisms on K3
surfaces}. Rev. Math. Iberoam. {\bf 29} (2013), 135--162.
%
\bibitem[vG01]{vG1} B.\ van Geemen, {\it Half twist of Hodge structures}, J. Math. Soc. Japan {\bf 53} (2001), 813--833.
%
\bibitem[vG92]{vG} B. van Geemen, {\it Projective models of Picard Modular Varieties},
Classification of Irregular Varieties, Springer LNM 1515 (1992) 68--99.
%
\bibitem[GH]{GH}
P. Griffiths, J. Harris \textit{Principles of algebraic geometry}.
John Wiley \& Sons Inc. New York, (1978).
%
\bibitem[Ha71]{Ha71}
W.J. Harvey, \textit{On branch loci in {T}eichm\"uller space}.
Trans. Amer. Math. Soc., \textbf{153}, (1971), 387--399.
%
\bibitem[K92]{Kondo} S. Kond$\bar{\rm o}$, {\it Automorphisms of algebraic K3 surfaces which act trivially on Picard groups}, J. Math. Soc. Japan, {\bf 44 }(1992), 75--98.
%
\bibitem[MSG]{MA} MAGMA Database of Small Groups;
http://magma.maths.usyd.edu.au/magma/htmlhelp/text404.htm.
%
\bibitem[M95] {M95}
R. Miranda, {\it Algebraic curves and {R}iemann surfaces}.
Graduate Studies in Mathematics, Vol {\bf 5}, American
Mathematical Society (1995).
%
\bibitem[M89] {M}
R.\ Miranda, {\it The Basic Theory of Elliptic Surfaces},  Dottorato di Ricerca in Matematica, Dipartimento di Matematica dell' Universit\`a di Pisa, ETS Editrice Pisa (1989). Available on line: http://www.math.colostate.edu/~miranda/BTES-Miranda.pdf.
%
\bibitem[MP10]{MP}
E. Mistretta, F. Polizzi, \textit{Standard isotrivial fibrations
with $p_g=q=1$ II}. J. Pure Appl. Algebra, \textbf{214},
(2010), 344--369.
%
\bibitem[P10] {P10}
F. Polizzi, \textit{Numerical properties of isotrivial
fibrations}. Geom. Dedicata,   \textbf{147}, (2010), 323--355.
%
\bibitem[OZ00]{OZ3}  K. Oguiso, D-Q. Zhang. {\it On VorontsovÕs theorem on K3 surfaces with non-symplectic group
actions}, Proc. Amer. Math. Soc. {\bf 128} (2000), 1571--1580.
%
\bibitem[OZ99]{OZ2} K. Oguiso, D-Q. Zhang, {\it K3 surfaces with order 11 automorphisms}, arXiv:math/9907020v1.
%
\bibitem[OZ98]{OZ1} K. Oguiso, D-Q. Zhang. {\it K3 surfaces with order five automorphisms}, J. Math. Kyoto Univ.
(1998) 419--438.
%
\bibitem[R79]{R} M. Reid, {\it Canonical 3-folds} Journ\'es de G\'eometrie Alg\'ebrique d'Angers, (1979),  273--310.
%
\bibitem[S96]{S96}
F. Serrano, \textit{Isotrivial fibred surfaces}. Annali di
Matematica pura e applicata,   \textbf{CLXXI}, (1996), 63--81.
%
\bibitem[SI77]{SI} T. Shioda, H. Inose,
{\it On singular K3 surfaces}, Complex analysis and algebraic geometry, Iwanami Shoten, Tokyo (1977) 119--136. 
\end{thebibliography}
\end{document}